\documentclass[review,3p]{elsarticle}

\usepackage{lineno,hyperref}
\usepackage{pgf}
\modulolinenumbers[5]

\journal{Journal of \LaTeX\ Templates}

\biboptions{sort&compress}

\usepackage{float}
\usepackage{lipsum}
\usepackage{amsfonts}
\usepackage{amsmath}
\usepackage{amsthm}
\usepackage{fourier}
\usepackage{xcolor}
\usepackage{natbib}
\usepackage{mathrsfs}
\usepackage{graphicx}
\usepackage{epstopdf}
\usepackage{algorithmic}
\usepackage[ntheorem]{empheq} 
\empheqset{box=\bigfbox}
\usepackage{stmaryrd}
\usepackage[final]{pdfpages}
\usepackage{subcaption}
\usepackage{verbatim}
\usepackage{booktabs}
\usepackage{amsopn}

\DeclareMathOperator*{\argmin}{arg\,min}

\newcommand{\Vn}{\mathbf{n}}

\newcommand{\Pol}{\mathbb{P}}

\newtheorem{ssmptn}{Assumption}

\newtheorem{prpstn}{Proposition}

\newtheorem{rmrk}{Remark}
\newtheorem{crllr}{Corollary}

\newcommand{\llbrace}{\lbrace \hspace{-.045cm}\lbrace }
\newcommand{\rrbrace}{\rbrace \hspace{-.05cm}\rbrace }

\newcommand{\lllvert}{\lvert \hspace{-.045cm}\lvert \hspace{-.045cm}\lvert}
\newcommand{\rrrvert}{\rvert \hspace{-.05cm}\rvert \hspace{-.05cm}\rvert}

\begin{document}

\begin{frontmatter}
\title{Incompressible flow modeling using an adaptive stabilized finite element method  based on residual minimization}


\author[six]{Felix E. Kyburg\corref{mycorrespondingauthor}}
\cortext[mycorrespondingauthor]{Corresponding author}
\ead{felixkyburg@gmail.com}

\author[one]{Sergio Rojas}
\author[one,three]{Victor M. Calo}

\address[one]{School of Earth and Planetary Sciences, Curtin University, Kent Street, Bentley, Perth, WA 6102, Australia}

\address[three]{Mineral Resources, Commonwealth Scientific and Industrial Research Organisation (CSIRO), Kensington, Perth, WA 6152, Australia}

\address[six]{Centro de Mec\'anica Computacional, Instituto Tecnol\'ogico de Buenos Aires, Madero, Buenos Aires, Argentina}

\begin{abstract}

We model incompressible flows with an adaptive stabilized finite element method Stokes flows, which solves a discretely stable saddle-point problem to approximate the velocity-pressure pair. Additionally, this saddle-point problem delivers a robust error estimator to guide mesh adaptivity. We analyze the accuracy of different discrete velocity-pressure pairs of continuous finite element spaces, which do not necessarily satisfy the discrete inf-sup condition. We validate the framework's performance with numerical examples.

\end{abstract}


\end{frontmatter}

\linenumbers

\section{Introduction}

The Stokes system describes the motion of incompressible fluid flows at low Reynolds numbers, where the effect of convection is negligible against the diffusion in the transport of momentum. The standard mixed finite element approach delivers a discrete saddle-point problem that requires compatible discrete spaces for the solution to satisfy the \textit{Babu\v{s}ka-Brezzi} inf-sup condition to achieve discrete stability~\cite{ Babuska1973, M2AN_1974__8_2_129_0}. Naive choices of Galerkin space pairs induce instabilities and locking~\cite{ ern2013theory}. In practice, richer mixed spaces like the mini-element~\cite{ Arnold1984} or the less efficient Taylor-Hood element are common. The mini-element enriches the velocity space to guarantee the inf-sup condition. Other strategies exist to circumvent this instability. Stabilized methods such as the Streamline upwind/Petrov-Galerkin (SUPG)~\cite{ BROOKS1982199, SUPG_stokes}, Galerkin/Least-Squares (GaLS)~\cite{ HUGHES1989173, 1989feaf.proc.1067F, HUGHES198785},  and the variational multiscale (VMS) finite element method~\cite{ LIU2006580} add residual-based terms to the classical weak discrete formulation to guarantee stability. Alternatively, the Discontinuous Petrov-Galerkin (DPG) method~\cite{ DEMKOWICZ20101558, demkowicz2011class, demkowicz2012class, zitelli2011class} exploits residual minimization to address this problem~\cite{ ROBERTS2014966}.  Alternatively, discontinuous Galerkin (DG) methods introduce discontinuities at element interfaces and use penalty terms to stabilize the solution~\cite{ di2011mathematical, doi:10.1137/1.9780898717440, kanschat_DG_flow}. For example, for the Stokes problem, \citet{ 10.2307/4100956} proposed equal-order methods for the velocity and pressure, while Girault et al.~\cite{ Wheeler_Stokes} proposed methods where the velocity has an extra order compared to the pressure order.

In this paper, we model Stokes flows using an adaptive stabilized finite element method~\cite{ SergiosMethod, 10.1007/978-3-030-50417-5_15} and study its performance on simplicial elements combined with adaptive mesh refinement.  The methodology inherits the inf-sup stability from the DG methods to construct a stable residual minimization. Within this framework, we can consider pairs of finite element spaces for the solution that are not compatible in the sense of the Babu\v{s}ka-Brezzi condition~\cite{ Babuska1973, M2AN_1974__8_2_129_0}. For example, equal-order velocity-pressure pairs deliver stable discrete solutions. Starting from a stable DG pair, we solve a residual minimization problem where the solution fields are continuous. The crucial insight is that we measure the error in a dual norm of the DG space, which results in an inf-sup stable formulation in terms of the DG norm. We can also use independent residual representations to measure the discrete errors of the continuous velocity and pressure approximations. Therefore, we can design adaptive refinement strategies that can take into account the pressure and velocity solutions separately.

The symmetric interior penalty method (SIPG)~\cite{ doi:10.1137/0719052, doi:10.1137/0715010} has optimal $L^2$ velocity convergence. Nevertheless, our trial spaces zero out several terms from the standard DG bilinear forms, which reduce the convergence suboptimal for even velocity polynomial degrees on uniform meshes  (see~\cite{  doi:10.1137/1.9780898717440}). Nevertheless, the refinement strategy we propose delivers optimal convergence rates for all polynomial orders.

In Sections~\ref{stokes_intro} and~\ref{dis_setting}, we present the Stokes problem and its discrete setting. In Section~\ref{stab_method}, we describe the new adaptive stabilized finite element framework for the Stokes equations in Section. In Section~\ref{num_aspects}, we briefly discuss some numerical aspects, while Section~\ref{num_examples} describe the numerical examples that show the impact of adaptivity. Finally, we detail our conclusions in Section~\ref{conclusions}.

\section{The Stokes Problem}\label{stokes_intro}
Let $\Omega$ be an open bounded domain ($\Omega \subset \mathbb{R}^d$ with $d=2,3$) with boundary $\Gamma = \partial \Omega$. The Stokes flow problem consists of finding a field function $\boldsymbol{u}$ (velocity), and a scalar function $p$ (pressure) satisfying:
\begin{subequations}\label{cont_stokes}
	\begin{align}
	-\nu \Delta \boldsymbol{u} + \nabla p &= \boldsymbol{f} \qquad & \text{in } \medspace  \Omega, \label{eq_momentum} \\ 
	\nabla \cdot \boldsymbol{u} &= 0 \qquad &\text{in } \medspace \Omega, \label{eq_conservation} \\ 
	\boldsymbol{u} &= \boldsymbol{u_0} \quad &\text{on } \medspace \Gamma. \label{eq_bcDirichlet}
	\end{align}
\end{subequations}

In the above equations, \eqref{eq_momentum} denotes the \textit{momentum balance}, while~\eqref{eq_conservation}, the \textit{mass balance} of the flow and \eqref{eq_bcDirichlet} the Dirichlet boundary conditions. We express the body force acting on the fluid as $\boldsymbol{f}$ and the kinematic viscosity of the medium as $\nu > 0$. Without loss of generality, we set~$\nu$ equal to~1. The pressure solution of  problem~\eqref{cont_stokes} is unique up to a constant. Therefore, to ensure uniqueness for~\eqref{cont_stokes}, we assume that $\langle p \rangle_{\Omega}=0$, where $\langle \cdot \rangle_{\Omega}$ denotes the mean value over $\Omega$: 
\begin{equation}\label{mean_pressure}
\langle p \rangle_{\Omega}= \frac{1}{|\Omega|}\int_{\partial\Omega}p.
\end{equation}

\subsection{Continuous weak formulation}

Denoting by $L^2(\Omega)$ the space of square integrable functions defined on $\Omega$, and by $H^{1}(\Omega) = \{ v \in L^2(\Omega) : \nabla v \in [L^2(\Omega)]^d\}$. We introduce the velocity space $U = [H^{1}_{0}(\Omega)]^d$ and the pressure space as~$P = L^2_0(\Omega)$. Where the  space $H^{1}_{0}(\Omega)$ is defined as:
\begin{equation}
H^{1}_{0}(\Omega) := \{v \in H^1(\Omega)\medspace|\medspace v_{|\partial\Omega} = 0\}.
\end{equation}
and taking into account the condition~\eqref{mean_pressure}, we define the space for the pressure as follows:
\begin{equation}
L^2_0(\Omega) := \{q \in L^2(\Omega)\medspace|\medspace \langle q \rangle_{\Omega} = 0\}.
\end{equation}
We introduce the Hilbert mixed space $\boldsymbol{X} := U \times P = [H^{1}_{0}(\Omega)]^d \times L^2_0(\Omega)$  doted with the following norm:
\begin{equation}
\parallel(\boldsymbol{v},q)\parallel_{\boldsymbol{X}} := \left(\parallel \mathbf{v} \parallel^2_U + \parallel q \parallel^2_P   \right)^{\frac{1}{2}}
\end{equation}
where $\| \cdot \|_U \text{ and } \|\cdot\|_P$ are defined as:
\begin{subequations}
	\begin{align}
	\parallel \boldsymbol{v}\parallel_U &:= \parallel v \parallel_{[H^1(\Omega)]^d} = \left(\sum^d_{i=1} \parallel v_i \parallel^2_{H^1(\Omega)} \right)^{\frac{1}{2}} = \left(\sum^d_{i=1} \left( \parallel v_i \parallel^2_{L^2(\Omega)} + \parallel \nabla v_i \parallel^2_{L^2(\Omega)}\right) \right)^{\frac{1}{2}} ,\\
	\parallel q\parallel_P &:= \parallel q \parallel_{L^2(\Omega)}.        
	\end{align}
\end{subequations}

The weak variational formulation for the Stokes problem~\eqref{cont_stokes} is well posed (c.f.~\cite{di2011mathematical}) and reads:
\begin{equation}\label{weakFormulation}
\left\{\begin{array}{l}
\text{Find } (\boldsymbol{u},p) \in \boldsymbol{X} := U \times P \text{, such that:} \smallskip \\
a(\boldsymbol{u},\boldsymbol{v}) + b(\boldsymbol{v},p) = l(\boldsymbol{v}), \qquad \forall \boldsymbol{v} \in U, \smallskip \\
d(\boldsymbol{u},q) = 0, \qquad  \forall q \in P.
\end{array}
\right.
\end{equation}
where the bilinears forms in~\eqref{weakFormulation}are defined as:
\begin{subequations}\label{bilinear_weakFormulation}
	\begin{align}
	a(\boldsymbol{v},\boldsymbol{w}) &:= \int_{\Omega}  \medspace \nabla \boldsymbol{v} : \nabla \boldsymbol{w} \medspace  = \sum_{i,j=1}^{d} \int_{\Omega} \partial_j v_i \partial_{j} w_i = \left(\nabla v, \nabla w\right)_{[L^2(\Omega)]^{d,d}},\\
	b(\boldsymbol{v},q)  &:= -\int_{\Omega} q \medspace \nabla \cdot \boldsymbol{v} \medspace ,\\
	d(\boldsymbol{v},q) &:= \int_{\Omega} \boldsymbol{v}  \medspace \cdot \nabla  q, \medspace 
	\end{align}
\end{subequations}
and the linear form:
\begin{equation}\label{linear_weakFormulation}
l(\boldsymbol{v})   := \int_{\Omega} \boldsymbol{f} \cdot \boldsymbol{v} \medspace  = \sum_{i=1}^{d} \int_{\Omega} f_i v_i.
\end{equation}

\section{Discrete setting}\label{dis_setting}
We partition the domain $\Omega$ into a conforming mesh $\Omega_h$ consisting of $N$ open disjoint elements $K_i$, satisfying
\begin{equation}
\Omega_h := \bigcup^N_{i=1} K_i \text{ satisfies } \text{int}(\overline{\Omega_h}) = \Omega.  
\end{equation}

For each element $K_i$, we denote by $\partial K_i$ its boundary. We denote a mesh face as $F$, and we assume it belonging to one, and only one, of the following two sets:
\begin{itemize}
\item {\em interior faces:} denoted by $F^i_h$, if $F = \partial K_i \cap \partial K_j$, where $K_i \text{ and } K_j$ are two mesh elements.
\item {\em boundary faces:} denoted by $F^b_h$, if $F = \partial K_i \cap \partial \Omega$.
\end{itemize}
Furthermore, we define the skeleton of $\Omega_h$ as $F_h := F^i_h \cup F^b_h$. 

\begin{figure}[h!]
	\centering
	\includegraphics[width=0.3\textwidth]{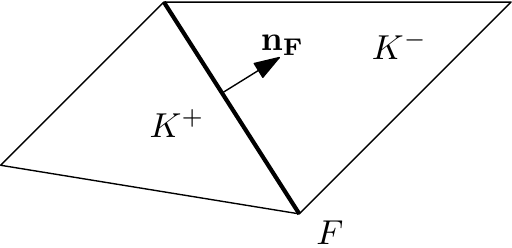}
	\caption{Skeleton orientation over an interface $F = \partial K^+ \cap \partial K^-$.}
	\label{fig:skeleton}
\end{figure}

For a given function $v$ sufficiently regular, we define the jump $\llbracket v \rrbracket_F$ and the average $\llbrace v \rrbrace_F$ operators, over an interior face $F$, as:
\begin{align}
\llbracket v \rrbracket_F(x) &:= v^+(x) - v^-(x),\label{eq:jump1}\\
\llbrace  v \rrbrace_F(x) &:= \dfrac{1}{2} (v^+(x) + v^-(x)), \label{eq:jump2}
\end{align}
where $v^{\pm}$ denotes the traces over $F$ following a predefined normal $\Vn_F$ over $F$ (see Figure~\ref{fig:skeleton}).

We extend the definitions \eqref{eq:jump1} and \eqref{eq:jump2} to $F \in F^b_h$, as $\llbracket v \rrbracket_F(x) := \llbrace  v \rrbrace_F(x) := v_F $.

We define the following piece-wise polynomial spaces defined over $\Omega_h$: 
\begin{subequations}\label{polynomial_def}
	\begin{align}
	\mathbb{P}^{k}(\Omega_h)    &:= \{v_h \in C^0(\overline{\Omega}): v_h|_{K_i} \in \mathbb{P}^k(K_i), \forall i = 1,...,N \}\label{polynomial_def1} \\
	\mathbb{P}^{k}_d(\Omega_h) &:= \{v_h \in L^2(\Omega): v_h|_{K_i} \in \mathbb{P}^k(K_i), \forall i = 1,...,N \}\label{polynomial_def_broken}.
	\end{align}
\end{subequations}

For a discrete space $W_h$, defined in terms of the discrete spaces~\eqref{polynomial_def1} and~\eqref{polynomial_def_broken}, doted with a discrete norm $\| \cdot \|_{W_h}$, we denote by $W_h^\ast$ the dual space of $W_h$, and we define the dual norm $\|\cdot\|_{W_h^\ast}$ as: 
\begin{equation}
\displaystyle \| \Psi_h \|_{W^{\ast}_h} := \sup_{0 \neq w_h \in W_h} \frac{\langle \Psi_h,w_h \rangle_{W^{*}_h \times W_h}}{\|w_h\|_{W_h}} \qquad\forall \, \Psi_h \in W_h^\ast,   
\end{equation}
where $\langle \cdot, \cdot \rangle_{W^{*}_h \times W_h}$ denotes the duality pairing in $W^{*}_h \times W_h$. Finally, we denote by $R_{W_h}$ the Riesz map:
\begin{equation}\label{eq:riesz_w}
  \begin{array}{rcl}
    R_{W_h} &:& W_h \rightarrow W_h^\ast \smallskip\\
            && \left< R_{W_h} y_h, w_h\right>_{W_h^\ast \times W_h} := (y_h,w_h)_{W_h} , \qquad \forall w_h \in W_h.
  \end{array}
\end{equation}
The Riesz map is an isometric isomorphism satisfying (see Theorem 6.4.1 in~\cite{  2358053}):
\begin{equation}\label{riesz_identity}
  \displaystyle \| \Psi_h z_h \|_{W^*_h} := \|R^{-1}_{W_h}(\Psi_h z_h)  \|_{W_h} \qquad \forall \Psi_h \in W_h,
\end{equation}
where $R^{-1}_{W_h} : W_h^\ast \rightarrow W_h$ denotes the inverse of the Riesz map~\eqref{eq:riesz_w}.

\subsection{Non-conforming discontinuous Galerkin formulation for the Stokes problem}\label{DGF}

In this section, we present a discontinuous Galerkin discretization of the steady Stokes system. We define $\boldsymbol{Y_h} := V_h \times Q_h = [\mathbb{P}^{k}_d(\Omega_h)]^d \times~\mathbb{P}^{k}_{d,0}(\Omega_h)$, where the subspace $\mathbb{P}^{k}_{d,0}(\Omega_h) \subset \mathbb{P}^{k}_d(\Omega_h)$ solely includes functions that satisfy $\langle p \rangle_{\Omega}=~0$. We consider the following DG formulation (c.f.~\cite{ di2011mathematical}): 
\begin{equation}\label{DG_problem}
\left\{\begin{array}{l l l}
\text{Find } (\boldsymbol{u_h},p_h) \in \boldsymbol{Y_h} := V_h \times Q_h, &\text{such that:} \smallskip \\
a_h(\boldsymbol{u_h},\boldsymbol{v_h}) + b_h(\boldsymbol{v_h},p_h) &= l_h(\boldsymbol{v_h}), & \forall \boldsymbol{v_h} \in V_h, \smallskip \\
d_h(\boldsymbol{u_h},q_h)  - s_h(q_h,p_h) &= 0, & \forall q_h \in Q_h,
\end{array}
\right.
\end{equation}
where the corresponding bilinear forms are defined as follows.
The bilinear $a_h(\cdot,\cdot)$ corresponds to the Symmetric interior penalization Galerkin (SIPG) (see~\cite{ di2011mathematical, doi:10.1137/1.9780898717440}):
\begin{align}\label{a_sip_bilinear}
a_h(\boldsymbol{v_{h}},\boldsymbol{w_{h}}) = \int_{\Omega}  \medspace \nabla_h \boldsymbol{v_{h}} : \nabla
_h \boldsymbol{w_{h}} \medspace + \sum_{F \in F_h} \frac{\eta}{h_F} \int_{F} \llbracket \boldsymbol{v_{h}} \rrbracket \cdot \llbracket \boldsymbol{w_{h}} \rrbracket 
- \sum_{F \in F_h} \int_{F} (\llbrace \nabla_h \boldsymbol{v_{h}} \rrbrace  \mathbf{n_F} \cdot \llbracket \boldsymbol{w_{h}} \rrbracket +  \llbrace \nabla_h \boldsymbol{w_{h}} \rrbrace  \mathbf{n_F} \cdot \llbracket \boldsymbol{v_{h}} \rrbracket ), 
\end{align}
where $h_F$ is the local length scale (i.e., diameter of $F$) and the user-defined real constant $\eta > 0$, which denotes the penalty parameter. The bilinears $b_h(\cdot,\cdot)$ and $d_h(\cdot,\cdot)$ correspond to the pressure-velocity coupling:
\begin{subequations}
	\begin{align}
	b_h(\boldsymbol{v_h},q_h)  = - \int_{\Omega} q_h  \nabla_h \cdot \boldsymbol{v_h} + \sum_{F \in F_h} \int_F \llbracket \boldsymbol{v_h} \rrbracket \cdot \mathbf{n_F} \llbrace q_h \rrbrace , \\
	d_h(\boldsymbol{v_h},q_h) = \int_{\Omega} \boldsymbol{v_h} \cdot \nabla_h q_h - \sum_{F \in F^i_h} \int_F \llbrace \boldsymbol{v_h} \rrbrace \cdot \mathbf{n_F} \llbracket q_h \rrbracket .  
	\end{align}  
\end{subequations}

The $s_h(\cdot,\cdot)$ term is a stabilization bilinear form controlling the pressure jumps. It is defined as:
\begin{equation}
s_h(q_h,r_h) =  \sum_{F \in F^i_h} h_F \int_F \llbracket q_h \rrbracket \llbracket r_h \rrbracket .
\end{equation}
This stabilization bilinear form $s_h(\cdot,\cdot)$ allows to consider spaces  with the same polynomial degree for both, velocities and pressures. Finally, the linear form $l_h(\boldsymbol{v_h})$ is
\begin{equation}\label{linear_DG}
l_h(\boldsymbol{v_h})   := \int_{\Omega} \boldsymbol{f} \cdot \boldsymbol{v} \medspace  -\int_{\partial\Omega} \boldsymbol{u_0} \cdot  (\nabla_h \boldsymbol{v_{h}} )  n_F \medspace + \sum_{F \in F^b_h} \frac{\eta}{h_F} \int_{F}  \boldsymbol{u_0}  \cdot  \boldsymbol{v_{h}} .
\end{equation}

We write problem~\eqref{DG_problem} in a simplified form as:
\begin{equation}\label{DG_problem2}
\left\{\begin{array}{l l l}
\text{Find } (\boldsymbol{u_h},p_h) \in \boldsymbol{Y_h} &\text{such that:} \smallskip \\
c_h((\boldsymbol{u_h},p_h),(\boldsymbol{v_h},q_h)) &= l(\boldsymbol{v_h}), & \forall (\boldsymbol{v_h},q_h) \in \boldsymbol{Y_h}, \smallskip \\
\end{array}
\right.
\end{equation}
where
\begin{equation}
c_h((\boldsymbol{u_h},p_h),(\boldsymbol{v_h},q_h)):= a_h(\boldsymbol{u_h},\boldsymbol{v_h}) + b_h(\boldsymbol{v_h},p_h) + d_h(\boldsymbol{u_h},q_h) - s_h(p_h,q_h).
\end{equation}

\begin{rmrk}[Pressure space for discontinuous Galerkin formulation]
		In this case, we consider the equal-order space $\boldsymbol{Y_h} := V_h \times Q_h = [\mathbb{P}^{k}_d(\Omega_h)]^d \times~\mathbb{P}^{k}_{d,0}(\Omega_h)$ for the discontinuous Galerkin formulation. The pressure subspace $\mathbb{P}^{k-1}_{d,0}(\Omega_h)$ can also be considered in the formulation, thus the $s_h(\cdot,\cdot)$ term is not needed. In that case, the same optimal convergence is obtained but the method is less efficient than the equal-order formulation (see \cite{10.2307/4100956}).
\end{rmrk}

\subsection{Discrete norms}

For the broken polynomial space $\boldsymbol{Y_h}:= V_h \times Q_h$, we consider the following norm (see~\cite{ di2011mathematical,doi:10.1137/1.9780898717440}):
\begin{equation}\label{eq:DG_stokes_norm}
\lllvert (\boldsymbol{v_h},q_h) \rrrvert := \left(\lllvert \boldsymbol{v_h}\rrrvert^2_{v} + \|q_h\|^2_P + |q_h|^2_p \right)^{\frac{1}{2}},
\end{equation}
with $\lllvert \boldsymbol{v_h}\rrrvert_{v}$ given by:
\begin{equation}
\lllvert \boldsymbol{v_h} \rrrvert_{v} :=   \left( \|\nabla_h v_h \|^2_{[L^2(\Omega)]^{d,d}} + |v_h|^2_J \right)^{\frac{1}{2}},
\end{equation}
%
and the $| \cdot |_J$-seminorm acting on vector valued arguments as:
\begin{equation}
|v_h|_J = \left( \sum_{F \in F_h} \frac{\eta}{h_F} \left\| \llbracket v_h \rrbracket  \right\|^2_{[L^2(F)]^d}   \right)^{\frac{1}{2}}.
\end{equation}
\noindent
For the pressure related terms, $\|q_h \|_P$ is
\begin{equation}
\|q_h\|_P := \|q_h\|_{[L^2(\Omega)]}, 
\end{equation}
and $|q_h|_p$ is the following pressure semi-norm
\begin{equation}
|q_h|_p := \left(\sum_{F \in F^i_h} h_F \left\| \llbracket q_h \rrbracket \right\|^2_{L^2(F)} \right)^{\frac{1}{2}}.
\end{equation}

Problem~\eqref{DG_problem2} satisfies the following (see Lemma 6.13 in~\cite{di2011mathematical}): 
\begin{prpstn}\label{as:inf-sup}
	(Inf-sup stability): If the parameter $\eta \geq 0$ in the SIP bilinear form~\eqref{a_sip_bilinear} is large enough, then there exists a mesh independent constant $C_{\emph{sta}}>0$ such that, for all $(\boldsymbol{v_h},q_h) \in \boldsymbol{Y_h}$, it holds:
	\begin{equation}\label{eq:infsup_h}
	\displaystyle \sup_{0\neq (\boldsymbol{v_h},q_h) \in \boldsymbol{Y_h}} \dfrac{c_h((\boldsymbol{v_h},q_h),(\boldsymbol{w_h},r_h))}{\lllvert (\boldsymbol{w_h},r_h) \rrrvert}  \geq C_{\emph{sta}}\lllvert (\boldsymbol{v_h},q_h) \rrrvert.
	\end{equation}
\end{prpstn}

\begin{prpstn} \label{as:reg-consi} (Consistency with regularity) If the exact solution $(\boldsymbol{u},p)$ of~\eqref{weakFormulation} belongs to the extended space $\boldsymbol{Y_\#} := V_{\#} \times Q_{\#} $, with
	\begin{equation}
	V_{\#} := U \cap [H^2(\Omega)]^d , \qquad Q_{\#} := P \cap [H^1(\Omega)],
	\end{equation}
	then, the discrete bilinear form $c_h$ extends continuously to $\boldsymbol{Y_{h,\#}} \times \boldsymbol{Y_h}$, with $\boldsymbol{Y_{h,\#}} := \boldsymbol{Y_{\#}} + \boldsymbol{Y_h}$. Moreover, it holds:
	\begin{equation}\label{eq:consistency}
	c_h((\boldsymbol{u},p),(\boldsymbol{v_h},q_h)) = l(\boldsymbol{v_h}), \quad \forall \, (\boldsymbol{v_h},q_h) \in \boldsymbol{Y_h}.
	\end{equation}
\end{prpstn}

\begin{prpstn} \label{as:bound}
	(Boundedness): There exists a constant $C_{\emph{bnd}}>0$, uniform with respect to the mesh size, such that:
	\begin{equation}\label{eq:continuity}
	c_h((\boldsymbol{w},r),(\boldsymbol{v_h},q_h)) \leq C_{\emph{bnd}} \, \lllvert(\boldsymbol{w},r)\rrrvert_{\#} \lllvert (\boldsymbol{v_h},q_h) \rrrvert, \quad \forall \, ((\boldsymbol{w},r), (\boldsymbol{v_h},q_h)) \in \boldsymbol{Y_{h,\#}} \times \boldsymbol{Y_{h}}.
	\end{equation}
	
\end{prpstn}

\begin{prpstn}
	(A priori DG error estimate)
	There exists a unique $(\boldsymbol{u_h},p_h) \in \boldsymbol{Y_h}$ solution to the discrete problem~\eqref{DG_problem} satisfying the following a priori estimate:
	\begin{equation}
	\lllvert(\boldsymbol{u-u_h},p-p_h)\rrrvert \leq C \inf_{(\boldsymbol{v_h},q_h) \in Y_h}\lllvert(\boldsymbol{u-v_h},p-q_h)\rrrvert_{\#}.
	\end{equation} 
	where $C$ is a constant independent of the mesh.
\end{prpstn}

\section{The conforming adaptive stabilized finite element method}\label{stab_method} 
\subsection{Residual minimization problem}

Using the stable formulation~\eqref{DG_problem} as a starting point. The AS-FEM method \cite{SergiosMethod} consists of considering a conforming (in $U \times P$) trial space $\boldsymbol{X_h} := U_h \times P_h \subset \boldsymbol{Y_h}$, and solving the following residual minimization problem:
\begin{equation}\label{eq:min_prob}
\left\{\begin{array}{r l}
\multicolumn{2}{l}{\text{Find } (\boldsymbol{u_h},p_h) \in \boldsymbol{X_h} \subset \boldsymbol{Y_h},  \text{ such that:} \smallskip} \\
\displaystyle \left(\boldsymbol{u_h},p_h\right)&= \displaystyle \argmin_{(\boldsymbol{w_h},r_h) \in \boldsymbol{X_h}} \frac{1}{2} \left\|l_h- A_h \, \boldsymbol{w_h} -B_h \, r_h \right\|^2_{V_h^\ast} + \frac{1}{2}\|-D_h \, \boldsymbol{w_h} \|^2_{Q_h^\ast} \\[10pt]
&\displaystyle \overset{\eqref{riesz_identity}}{=} \argmin_{(\boldsymbol{w_h},r_h) \in \boldsymbol{X_h}} \frac{1}{2}\|R^{-1}_{V_h}(l_h- A_h \boldsymbol{w_h} - B_h r_h)\|^2_{V_h}   + \frac{1}{2}\|R^{-1}_{Q_h}(-D_h \boldsymbol{w_h} )\|^2_{Q_h},
\end{array}
\right.
\end{equation}
where the discrete operators in~\eqref{eq:min_prob} are defined as:
\begin{equation}
\begin{array}{lll}
A_h : U_{h,\#} \rightarrow V_h^{\ast},&  \langle A_h \boldsymbol{w},\boldsymbol{v_h} \rangle_{V_h^{\ast} \times V_h} &:= a_h(\boldsymbol{w},\boldsymbol{v_h}), \\
B_h : P_{h,\#} \rightarrow V_h^{\ast},&   \langle B_h r, \boldsymbol{v_h} \rangle_{V_h^{\ast} \times V_h} &:= b_h(r,\boldsymbol{v_h}), \\
D_h : U_{h,\#} \rightarrow Q_h^{\ast},&   \langle D_h \boldsymbol{w}, q_h \rangle_{Q_h^{\ast} \times Q_h} &:= d_h(\boldsymbol{w},q_h),\\
A^{\ast}_h : V_{h,\#} \rightarrow U_h^{\ast},&  \langle A^{\ast}_h \boldsymbol{v},\boldsymbol{z_h} \rangle_{U_h^{\ast} \times U_h} &:= a_h(\boldsymbol{z_h},\boldsymbol{v}), \\
D^{\ast}_h : Q_{h,\#} \rightarrow U_h^{\ast},&   \langle D^{\ast}_h q, \boldsymbol{z_h} \rangle_{U_h^{\ast} \times U_h} &:= d_h(\boldsymbol{z_h},q), \\
B^{\ast}_h : V_{h,\#} \rightarrow P_h^{\ast},&   \langle B^{\ast}_h \boldsymbol{v}, r_h \rangle_{P_h^{\ast} \times P_h} &:= b_h(r_h,\boldsymbol{v}),
\end{array}
\end{equation}

\noindent
Following~\cite{ SergiosMethod}, defining the residual representatives:
\begin{subequations}\label{eq:e_h}
	\begin{align}
	\boldsymbol{e^u_h} &:= R^{-1}_{V_h} (l_h-A_h \boldsymbol{u_h} - B_h p_h) \in V_h, \\
	e^p_h &:= R^{-1}_{Q_h} (-D_h \boldsymbol{u_h}) \in Q_h, 
	\end{align}
\end{subequations}
to solve the problem~\eqref{eq:min_prob} is equivalent to find the quartet $(\boldsymbol{e^u_h},e^p_h, \boldsymbol{u_h},p_h) \in (V_h \times Q_h \times U_h \times P_h$)
\begin{subequations}
	\label{eq:mix_form}
	\begin{empheq}[left=\left\{,right=\right.,box=]{alignat=3}
	\label{eq:mix_forma1}
	&\,\, (\boldsymbol{e^u_h} \, , \, \boldsymbol{v_h})_{V_h} + a_h(\boldsymbol{u_h} \, , \, \boldsymbol{v_h}) + b_h(p_h \, , \, \boldsymbol{v_h}) && =  l_h(\boldsymbol{v_h}),   &\quad&\forall\, \boldsymbol{v_h} \in V_h, \\
	\label{eq:mix_forma2}
	&\,\, (e^p_h \, , \, q_h)_{Q_h} + d_h(\boldsymbol{u_h} \, , \, q_h) && =  0,   &\quad&\forall\, q_h \in Q_h, \\
	\label{eq:mix_formb1}
	&\,\,a_h(\boldsymbol{z_h} \, , \, \boldsymbol{e^u_h)} + b_h(\boldsymbol{z_h} \, , \, e^p_h) &&=  0,  &\quad&\forall\, \boldsymbol{z_h} \in U_h, \\
	\label{eq:mix_formb2}
	&\,\,d_h(r_h \, , \, \boldsymbol{e^u_h}) &&=  0,  &\quad&\forall\, r_h \in P_h.
	\end{empheq}
\end{subequations}

We simplify the problem~\eqref{eq:mix_form} as:
\begin{subequations}
	\label{eq:mix_formb}
	\begin{empheq}[left=\left\{,right=\right.,box=]{alignat=3}
	&\displaystyle \,\, g_h((\boldsymbol{e^u_h},e^p_h ) \, , \, (\boldsymbol{v_h},q_h)) + n_h((\boldsymbol{u_h},p_h ) \, , \, (\boldsymbol{v_h},q_h)) && =  l_h(\boldsymbol{v_h}),   &\quad&\displaystyle \forall\, (\boldsymbol{v_h},q_h) \in \boldsymbol{Y_h}, \\
	&\displaystyle n^{\ast}_h((\boldsymbol{e^u_h},e^p_h ) \, , \, (\boldsymbol{z_h},r_h)) &&=  0,  &\quad&\displaystyle \forall\, (\boldsymbol{z_h},r_h) \in \boldsymbol{X_h}.
	\end{empheq}
\end{subequations}
where
\begin{subequations}
	\begin{align}
	g_h((\boldsymbol{e^u_h},e^p_h )&:= (\boldsymbol{e^u_h} \, , \, \boldsymbol{v_h})_{V_h} + (e^p_h \, , \, q_h)_{Q_h} \\
	n_h((\boldsymbol{u_h},p_h )	&:=	 a_h(\boldsymbol{u_h} \, , \, \boldsymbol{v_h}) + b_h(p_h \, , \, \boldsymbol{v_h}) + d_h(\boldsymbol{u_h} \, , \, q_h) \\
	n^{\ast}_h((\boldsymbol{e^u_h},e^p_h ) &:= \,\,a_h(\boldsymbol{z_h} \, , \, \boldsymbol{e^u_h)} + b_h(\boldsymbol{z_h} \, , \, e^p_h) + \,\,d_h(r_h \, , \, \boldsymbol{e^u_h})
	\end{align}  
\end{subequations}

The advantage of solving problem~\eqref{eq:mix_form} with $\boldsymbol{X_h}  \subset \boldsymbol{Y_h}$ is the fact that not only the solution $(\boldsymbol{u},p)$ is found. At the same time, we obtain the residual representation functions $(\boldsymbol{e^u_h},e^p_h)$ solution, which deliver two useful independent error representations for $\boldsymbol{u} \text{ and } p$  and can be conveniently used as markers for an adaptive refinement scheme. The following results hold (see~\cite{ SergiosMethod}):

\begin{prpstn}[A priori bounds and error estimates]\label{th:FEMwDG}
	The solution to the mixed problem~\eqref{eq:mix_form} is unique. Furthermore, the following a priori bounds are satisfied:
	\begin{equation}\label{eq:bounds}
	\lllvert(\boldsymbol{e^u_h},e^p_h)\rrrvert\leq \lllvert l_h\rrrvert_{\ast} \qquad \hbox{ and }\qquad \lllvert(\boldsymbol{u_h},p_h) \rrrvert_{St} \leq \dfrac{1}{C_{\emph{sta}}}\lllvert l_h\rrrvert_{\ast}\,\,,
	\end{equation}
	Remembering that $(\boldsymbol{u},p) \in \boldsymbol{X_{\#}}$ is the exact solution to the  continuous problem~\eqref{weakFormulation}, the following a priori error estimate holds:
	\begin{equation}\label{eq:apriori}
	\displaystyle \lllvert(\boldsymbol{u}-\boldsymbol{u_h},p-p_h)\rrrvert \leq C \inf_{(\boldsymbol{w_h},r_h) \in \boldsymbol{X_h}}\lllvert(\boldsymbol{u}-\boldsymbol{w_h},p-r_h)\rrrvert_{{\#}}\,\,,
	\end{equation}
\end{prpstn}

\begin{prpstn}[Error representative] \label{th:err_rep}
	Let $(\boldsymbol{u_h},p_h) \in \boldsymbol{X_h}$ be the third and fourth components of the quartet  $(\boldsymbol{e^u_h},e^p_h, \boldsymbol{u_h},p_h) \in (V_h \times Q_h \times U_h \times P_h)$ solving~\eqref{eq:mix_form}. Let $(\boldsymbol{w_h},r_h) \in \boldsymbol{Y_h}$ be the unique solution to~\eqref{DG_problem}. Then, the following holds:
	\begin{equation} \label{eq:bnd_uh-thetah}
	\lllvert(\boldsymbol{u_h-w_h},p_h-r_h)\rrrvert \le \dfrac{1}{C_{\emph{sta}}} \lllvert(\boldsymbol{e^u_h},e^p_h)\rrrvert.
	\end{equation} 
\end{prpstn}

\begin{crllr}[Reliability]	
	Under the same hypothesis of Proposition~\ref{th:err_rep}, the following holds:
	\begin{equation} \label{eq:efficiency}
	\lllvert(\boldsymbol{e^u_h},e^p_h)\rrrvert  \le C_{\emph{bnd}}\lllvert(\boldsymbol{u-u_h},p-p_h)\rrrvert_{\#} .
	\end{equation} 
\end{crllr}

\begin{ssmptn}\label{saturation}
	There exists a real number $\delta\in [0,1)$, uniform with respect to the mesh size, such that:
	\begin{equation}
	\lllvert(\boldsymbol{u-w_h},p-r_h)\rrrvert \le \delta \lllvert(\boldsymbol{u-u_h},p-p_h)\rrrvert.
	\end{equation}
\end{ssmptn}

\begin{crllr}[Efficiency]	
	Under the same hypothesis of Proposition~\ref{th:err_rep}, if in addition Assumption~\ref{saturation} holds true, then:
	\begin{equation} \label{eq:a_posteriori}
	\lllvert(\boldsymbol{u-u_h},p-p_h)\rrrvert \le \dfrac{1}{(1-\delta)C_{\emph{sta}}} \lllvert(\boldsymbol{e^u_h},e^p_h)\rrrvert.
	\end{equation} 
\end{crllr}

\section{Numerical aspects}\label{num_aspects}
Implementing the system of equations of problem~\eqref{eq:mix_form} results in the following block matrix system:
\begin{equation}
\left[
\begin{array}{c c | c c}
G_u & 0 & A & B \\
0   & G_p & D & 0 \\
\hline
A^* & D^* & 0 & 0 \\
B^* & 0 & 0 & 0
\end{array}
\right] \cdot 
\left[
\begin{array}{c}
\boldsymbol{e^u} \\ \boldsymbol{e^p} \\ \boldsymbol{u} \\ \boldsymbol{p}
\end{array}
\right] =
\left[
\begin{array}{c}
\boldsymbol{l} \\ 0 \\ 0 \\ 0
\end{array}
\right],
\end{equation}
or to simplify the system, like problem~\eqref{eq:mix_formb}, in the following form:
\begin{equation}\label{block_systemB}
\left[
\begin{array}{cc}
G & N \\
N^* & 0 \\
\end{array}
\right] \cdot 
\left[
\begin{array}{c}
\boldsymbol{e} \\ \boldsymbol{s}
\end{array}
\right] =
\left[
\begin{array}{c}
\boldsymbol{L} \\ 0 
\end{array}
\right], 
\end{equation}
where $G$ is the matrix induced by the Discontinuous Galerkin norms, $N$ is related to the original Stokes problem, and $N^*$ involves the orthogonality constraints of the error estimates to the solution space. The vector $\boldsymbol{e}$ includes the error estimators of both, velocity and pressure, while the vector $\boldsymbol{s}$ includes $(\boldsymbol{u},p)$.

We follow the framework employed in~\cite{ Siefert06} and consider a fixed-point method to solve the system~\eqref{block_systemB}. This yields into the following iterative system:
\begin{subequations}
	\label{eq:fp_iter}
	\begin{empheq}[left=\left\{,right=\right.,box=]{alignat=3}
	\label{eq:fp_iter_u}
	\boldsymbol{s_{i+1}} &= S^{-1}N^* \left(F^{-1}E\boldsymbol{e_i} + F^{-1}\boldsymbol{L}\right)\\
	\label{eq:fp_iter_p}
	\boldsymbol{e_{i+1}} &= \left(F^{-1}E\boldsymbol{e_i} + F^{-1}\boldsymbol{L}\right) - F^{-1}N\boldsymbol{s_{i+1}} 
	\end{empheq}
\end{subequations}

By performing the splitting of $G = F - E$, we obtain $F$, where the matrix $F$ is a preconditioner. We approach the exact Schur complement of the system~\eqref{block_systemB} with $S =N^* F^{-1} N$.  For this case, we use the CHOLMOD library~\cite{ Chen:2008:ACS:1391989.1391995} and select $F$ to be the sparse Cholesky factorization.

To solve $\boldsymbol{s_{i+1}}$, we apply a preconditioner $P_{Schur}$ for the Schur complement constructed as follows:
\begin{equation}\label{Preconditioner}
P_{Schur} = \left[
\begin{array}{c  c}
K_{AMG} &  0\\

0 & Q_{M_p} \\
\end{array}
\right]. 
\end{equation}

We approximate the matrices $K_{AMG} \text{ and } Q_{M_p}$ by applying the algebraic multigrid method to the standard Galerkin matrix for the diffusion problem using the PyAMG library~\cite{ OlSc2018}, and taking the diagonal of the mass matrix for the pressure, respectively~\cite{ elman_silvester_wathen_2005}.

\begin{rmrk}[Pressure mean-value condition]	
	
	The zero mean-value condition of the pressure field ~\eqref{mean_pressure} is settled using Lagrange multipliers and introducing as an additional unknown the constant of the pressure field (see DOLFIN Python Demos~\cite{LoggWellsEtAl2012a}). 
\end{rmrk}

\section{Numerical examples}\label{num_examples}

We apply our stabilized method to three different test cases to study the robustness of the error estimator.  We consider an analytical solution to measure the solution quality for both $\boldsymbol{u}$ and $p$ simultaneously under uniform refinements. We then apply the adaptive refinement strategy to model a circular segment with a reentrant corner. Lastly, we analyse the lid-driven cavity problem with uniform and adaptive refinement.

\begin{figure}[h]
	\centering
	\includegraphics[width=0.5\textwidth]{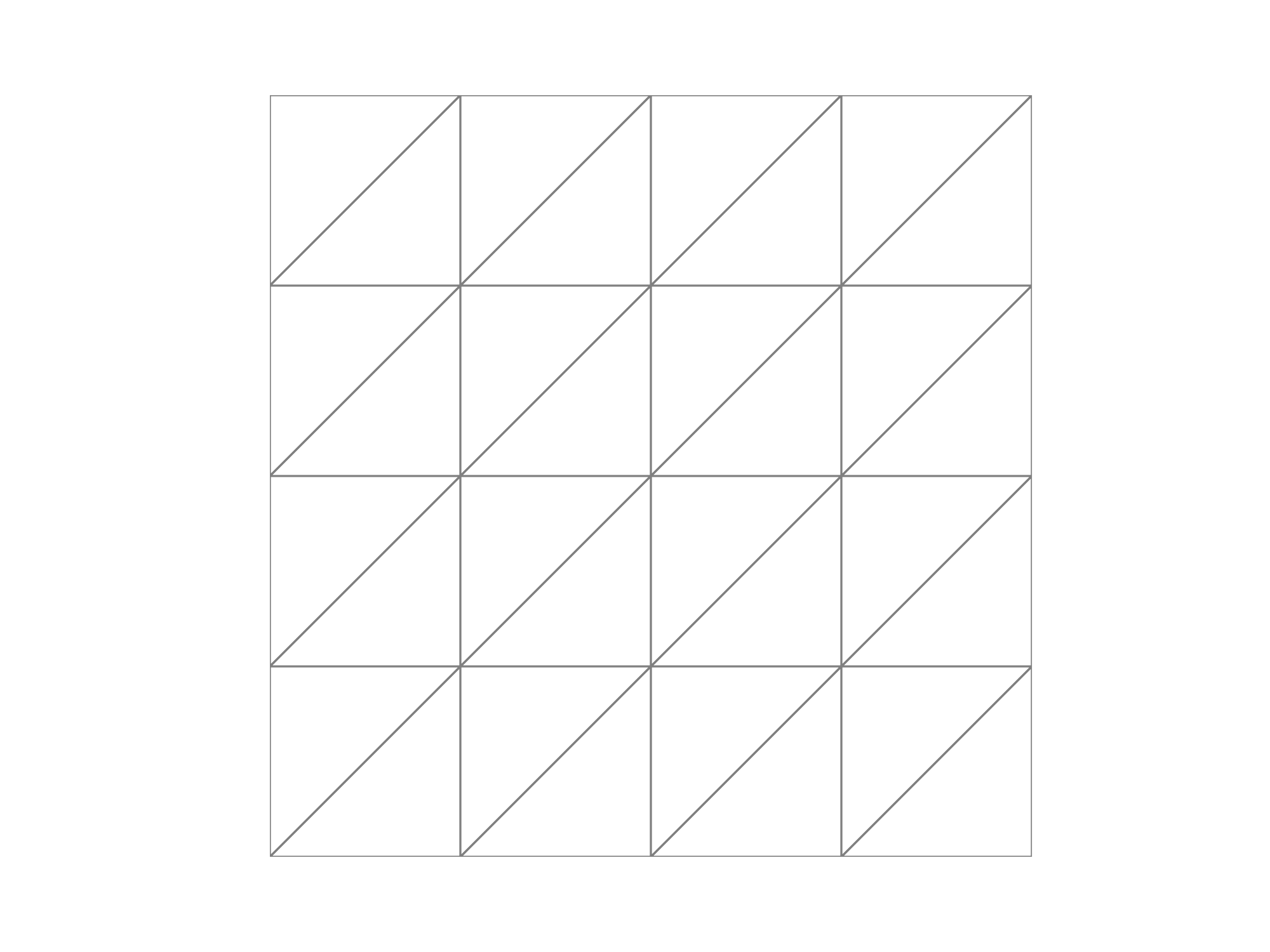}
	\caption{Initial mesh}
	\label{fig:ini_mesh1}
\end{figure}

For the analytic solution and lid-driven cases, we use a unit-square domain $\Omega = (0,1)^2$ and get an initial solution using an uniform triangular mesh (e.g. Figure~\ref{fig:ini_mesh1}). For the circular segment case, we use as initial mesh the circular segment shown in Figure~\ref{fig:case2_mesh0}. Next, in each subsequent level, we refine the mesh using a bisection-type refinement~\cite{ Bank1983SomeRA} and use the previous level solution as an initial guess. For the different cases, we use the computing platform FENICS~\cite{ AlnaesBlechta2015a} to solve the problem.

\subsection{Discretization setup}

We consider the broken polynomial spaces $V_h := [\Pol^{k}_{d}(\Omega_h)]^d$ and $Q_h := \Pol^{k}_{d}(\Omega_h)$, with $k=1,...,4$, defined in~\eqref{polynomial_def_broken} as test spaces. For the trial spaces, we consider several combinations of the polynomial spaces defined in~\eqref{polynomial_def} for $U_h$ and $P_h$. In particular, we label the trial space combinations of $U_h\times P_h:= [\mathbb{P}^k(\Omega_h)]^d \times \mathbb{P}_{0}^r(\Omega_h)$, used for the minimization problem~\eqref{eq:min_prob} as "$PkPr$," where~$k$ and~$r$ represent the polynomial order for the velocity and pressure, respectively. For the case $r=0$, the pressure space $P_h$ is $\mathbb{P}_{d,0}^0(\Omega_h)$.
In addition to the previous combinations, we label the original discontinuous Galerkin problem~\eqref{DG_problem} as "$DGk$," where $k$ is the polynomial order for the equal-order discontinuous velocity and pressure case. The spaces $V_h$  and $Q_h$ have the norms defined in~\eqref{eq:DG_stokes_norm}. 
%

%

\subsection{Adaptive refinement criteria}

We use an adaptive refinement criterion~\cite{SergiosMethod} that exploits the error representations $\boldsymbol{e^u_h} \in V_h$ and $e^p_h\in Q_h$ to drive the refinement process, which consists of four different steps in each iteration step:
\begin{enumerate}
	\item Solve the mixed problem~\eqref{eq:mix_form};
	\item Localize error over each element $K$ as $E_K := \rrrvert (\boldsymbol{e^u_h},e^p_h) \rrrvert_{K}$ with $\parallel \cdot \parallel_{V_h \times Q_h} = \rrrvert \cdot \rrrvert$. 

	\item We mark the elements to refine using the D\"{o}rfler bulk-chasing criterion~\cite{ doi:10.1137/0733054}. We mark the elements that contribute to the cumulative sum of the local value $E_K$ in decreasing order, while this sum remains below a user-specified fraction of the total value $\rrrvert(\boldsymbol{e^u_h},e^p_h)\rrrvert$. 
	\item To create the subsequent mesh, we refine the marked elements using a bisection-type refinement.
\end{enumerate}

\begin{figure}[ht!]
	\begin{subfigure}{0.5\textwidth}
		\includegraphics[width=\linewidth]{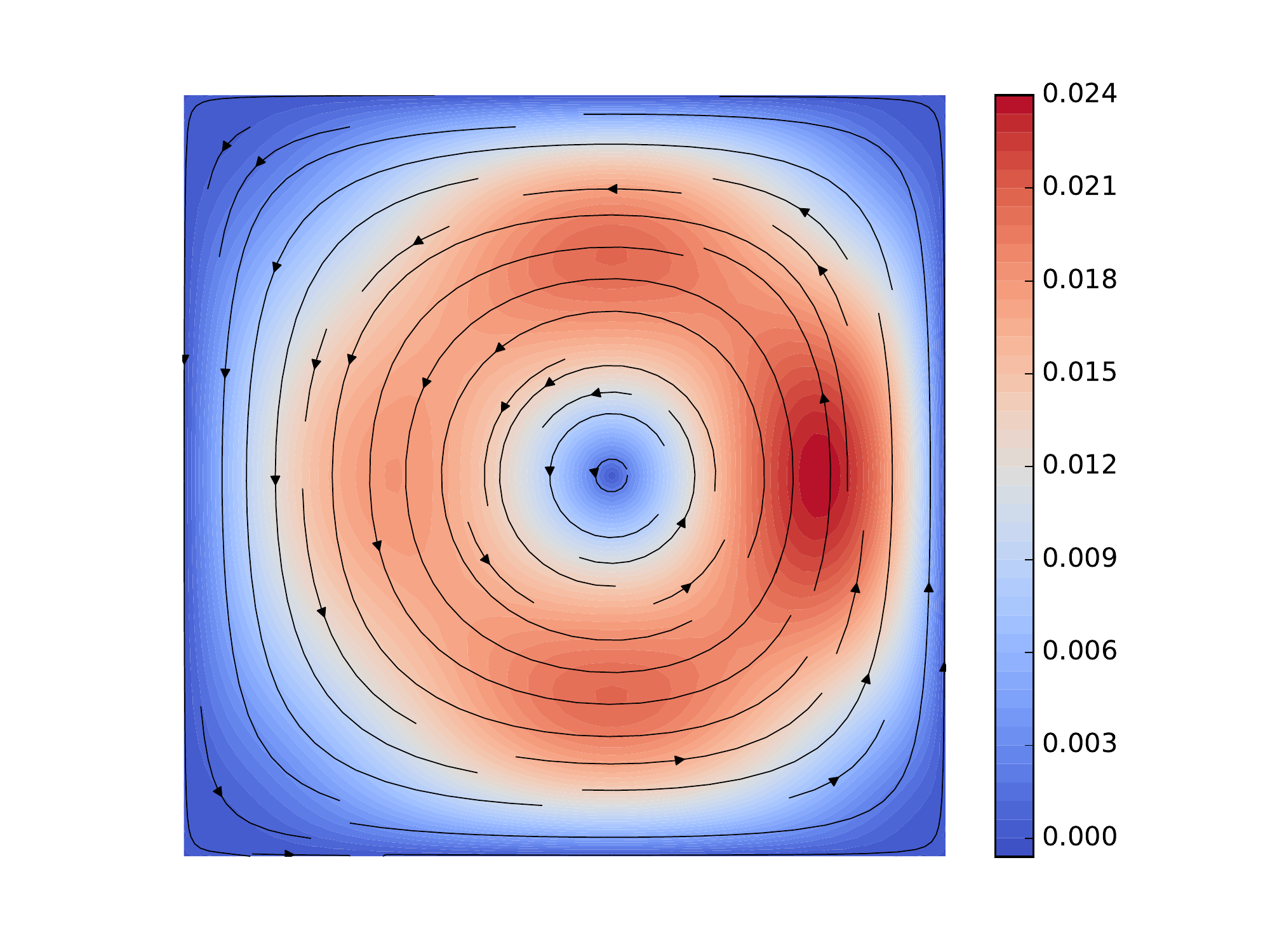}
		\caption{Velocity field ($\mathbf{u}$)}
		\label{fig:ref_case1_u}
	\end{subfigure}%
	\begin{subfigure}{0.5\textwidth}
		\centering
		\includegraphics[width=\linewidth]{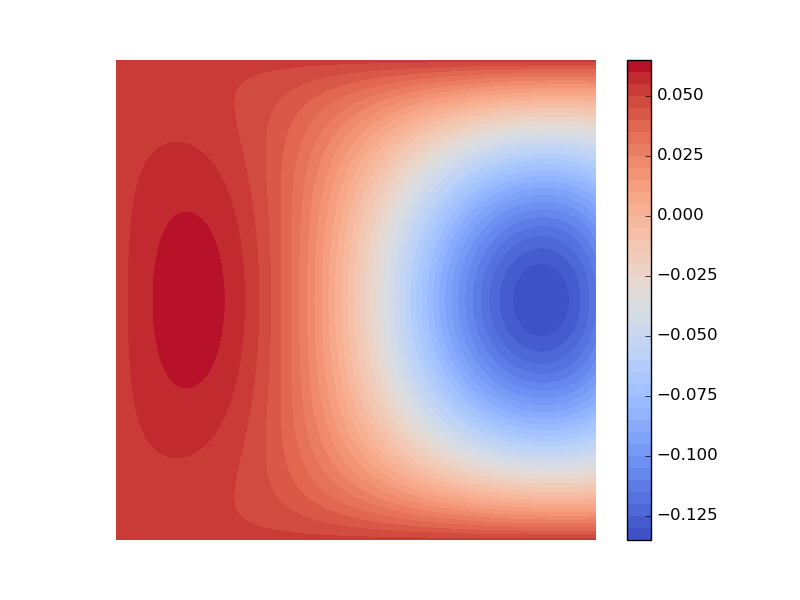}
		\caption{Pressre field ($p$)}
		\label{fig:ref_case1_p}
	\end{subfigure}
	\caption{Reference solution for the velocity and pressure fields in case $1$.}
	\label{fig:ref_case1}
\end{figure} 

\subsection{Convergence rates for uniform refinements for smooth
  analytic solution}\label{norm_selection}
  
We study a manufactured solution (c.f.~\cite{ doi:10.1002/fld.2337, CORTES2015123}) to validate and verify our implementation. The analytical expressions for the solution are:
%
\begin{equation}\label{u_exact_1}
\boldsymbol{u}_{sol} = 
\begin{bmatrix}
(2e^{x}(-1+x)^2 x^2 (y^2-y)(-1+2y)) \\
(-e^{x}(-1+x)x(-2+x(3+x))(-1+y)^2y^2)
\end{bmatrix},   
\end{equation}
\begin{equation}\label{p_exact_1}
\begin{split}
p_{sol} =& (-424+156e+(y^2-y)(-456+e^x(456+x^2(228-5(y^2-y)) \\
& + 2x(-228+(y^2-y)) + 2x^3(-36+(y^2-y)) + x^4(12+(y^2-y))))).
\end{split}
\end{equation}
Figure~\ref{fig:ref_case1} shows the streamlines, as well as the color plots of the solution.  Using the above velocity and pressure fields, we construct the body force term:
\begin{equation}
\boldsymbol{f}= -\Delta \boldsymbol{u}_{sol} + \nabla p_{sol}. 
\end{equation}

\begin{figure}[H]
  \centering
  \begin{subfigure}{0.45\textwidth}
    \centering
    \includegraphics[width=\linewidth]{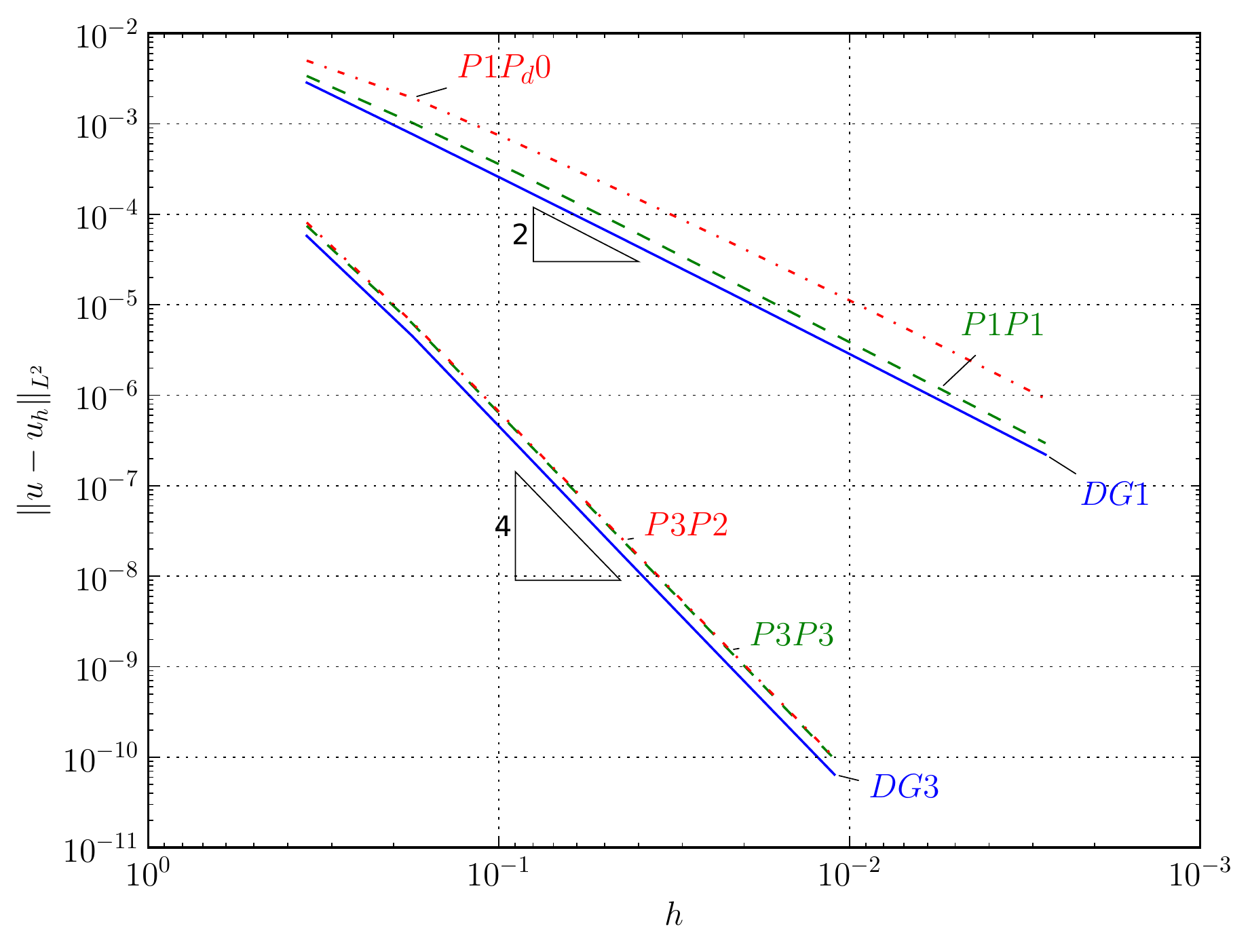}
    \caption{$k=1,3$}
    \label{fig:ref_case1_u_odd}
  \end{subfigure}
  \begin{subfigure}{0.45\textwidth}
    \centering
    \includegraphics[width=\linewidth]{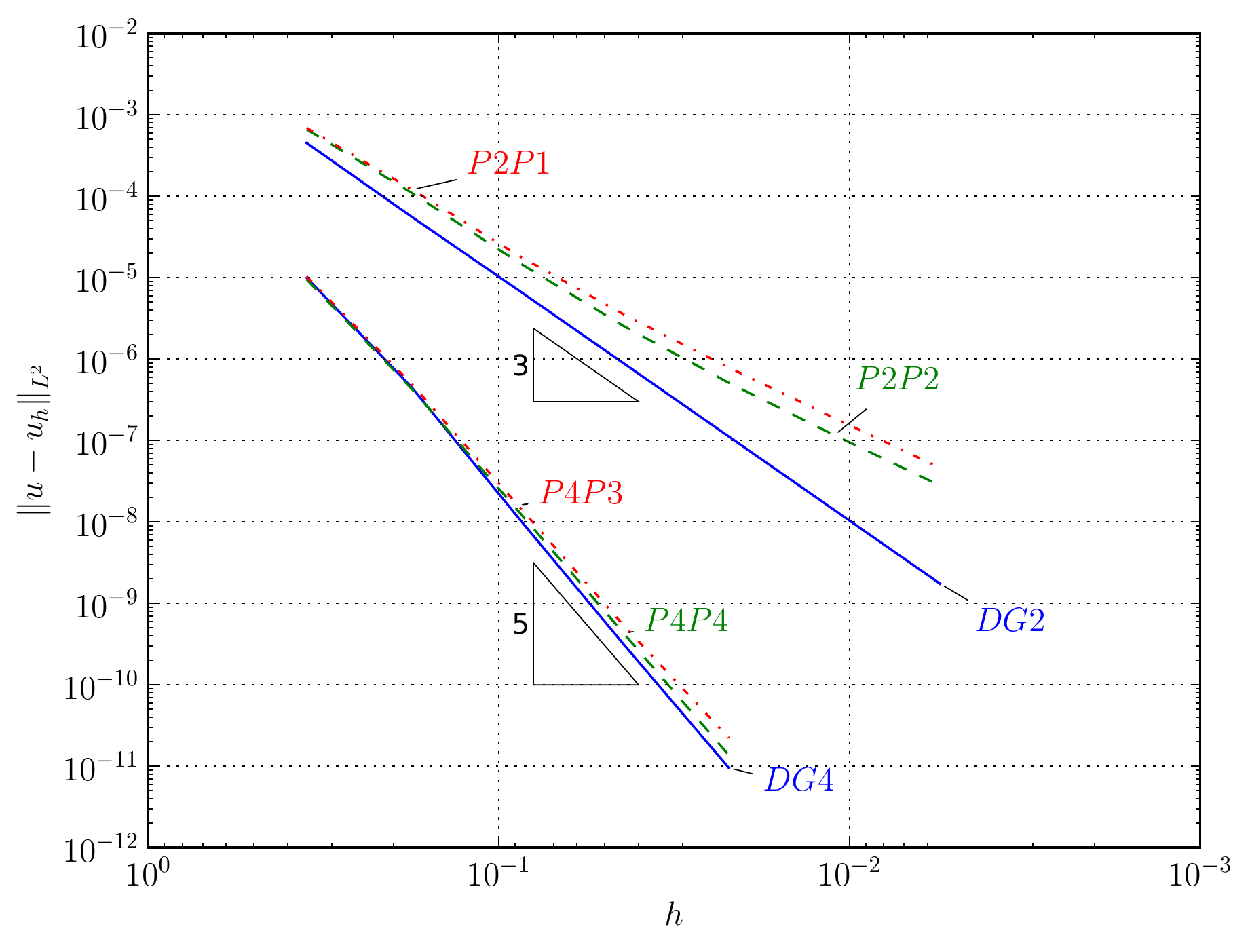}
    \caption{$k=2,4$}
    \label{fig:ref_case1_u_even}
  \end{subfigure}
  \caption{$\| \boldsymbol{u-u_h}\|_{L^2}$ for $k$ order of $\boldsymbol{u}$ and uniform refinement for a smooth solution.}
  \label{fig:case1_results_L2u}
\end{figure}

\begin{figure}[H]
  \centering
  \begin{subfigure}{0.45\textwidth}
    \centering
    \includegraphics[width=\linewidth]{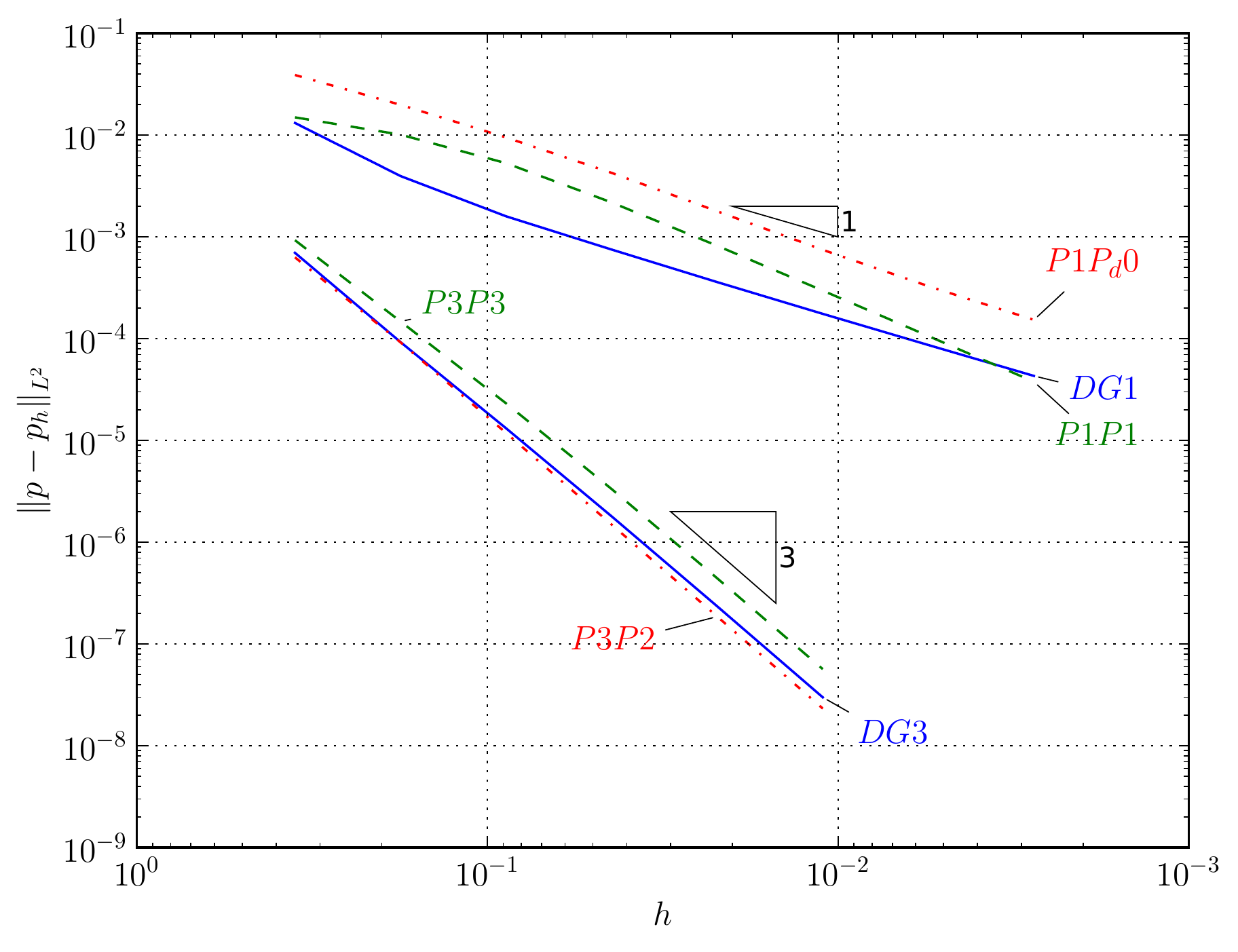}
    \caption{$k=1,3$}
    \label{fig:ref_case1_p_odd}
  \end{subfigure}
  \begin{subfigure}{0.45\textwidth}
    \centering
    \includegraphics[width=\linewidth]{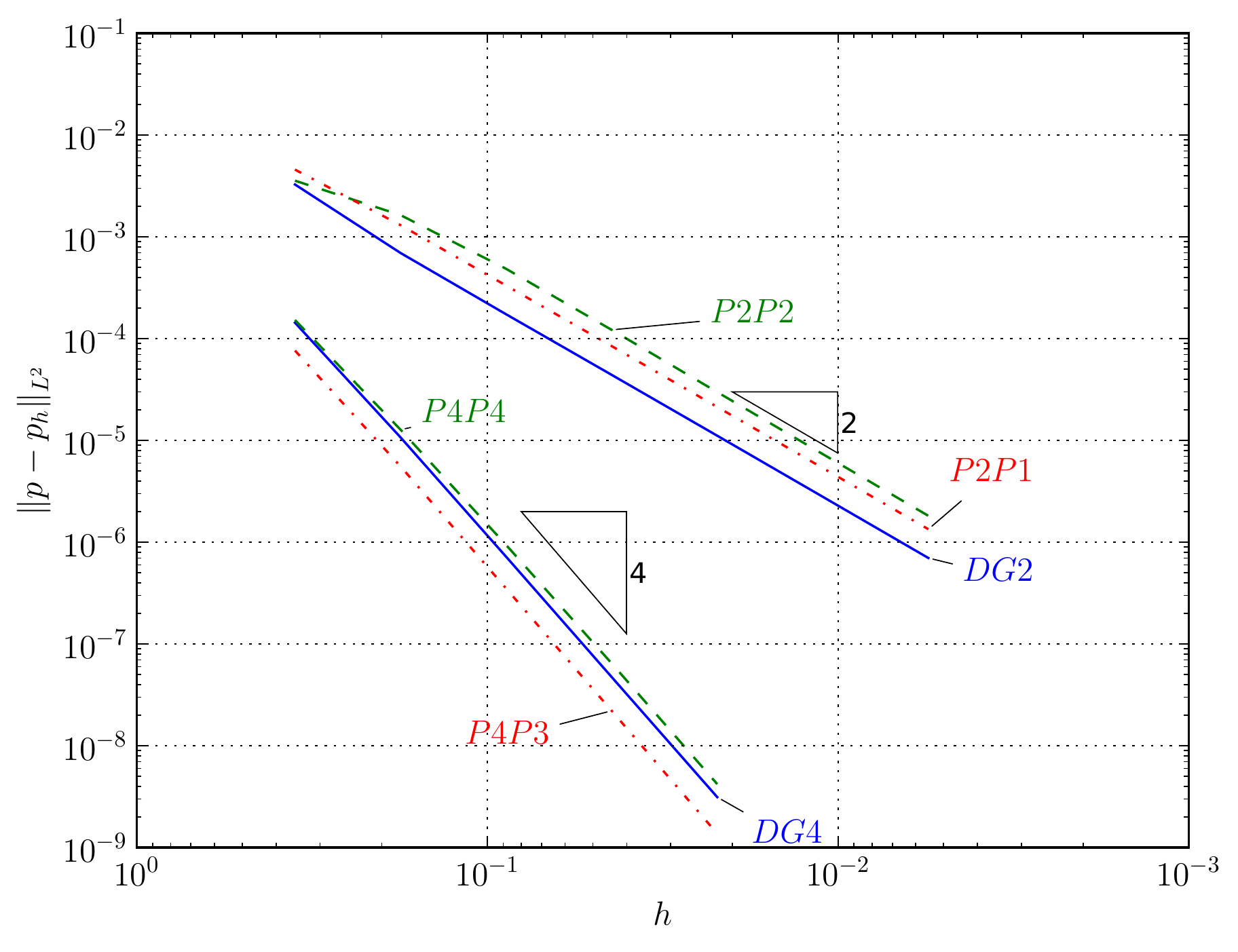}
    \caption{$k=2,4$}
    \label{fig:ref_case1_p_even}
  \end{subfigure}
  \caption{$\| p-p_h\|_{L^2}$ for $p$ for $k$ order of $\boldsymbol{u}$ and uniform refinement for a smooth solution.}
  \label{fig:case1_results_L2p}
\end{figure}

\begin{figure}[H]
  \centering
  \begin{subfigure}{0.45\textwidth}
    \centering
    \includegraphics[width=\linewidth]{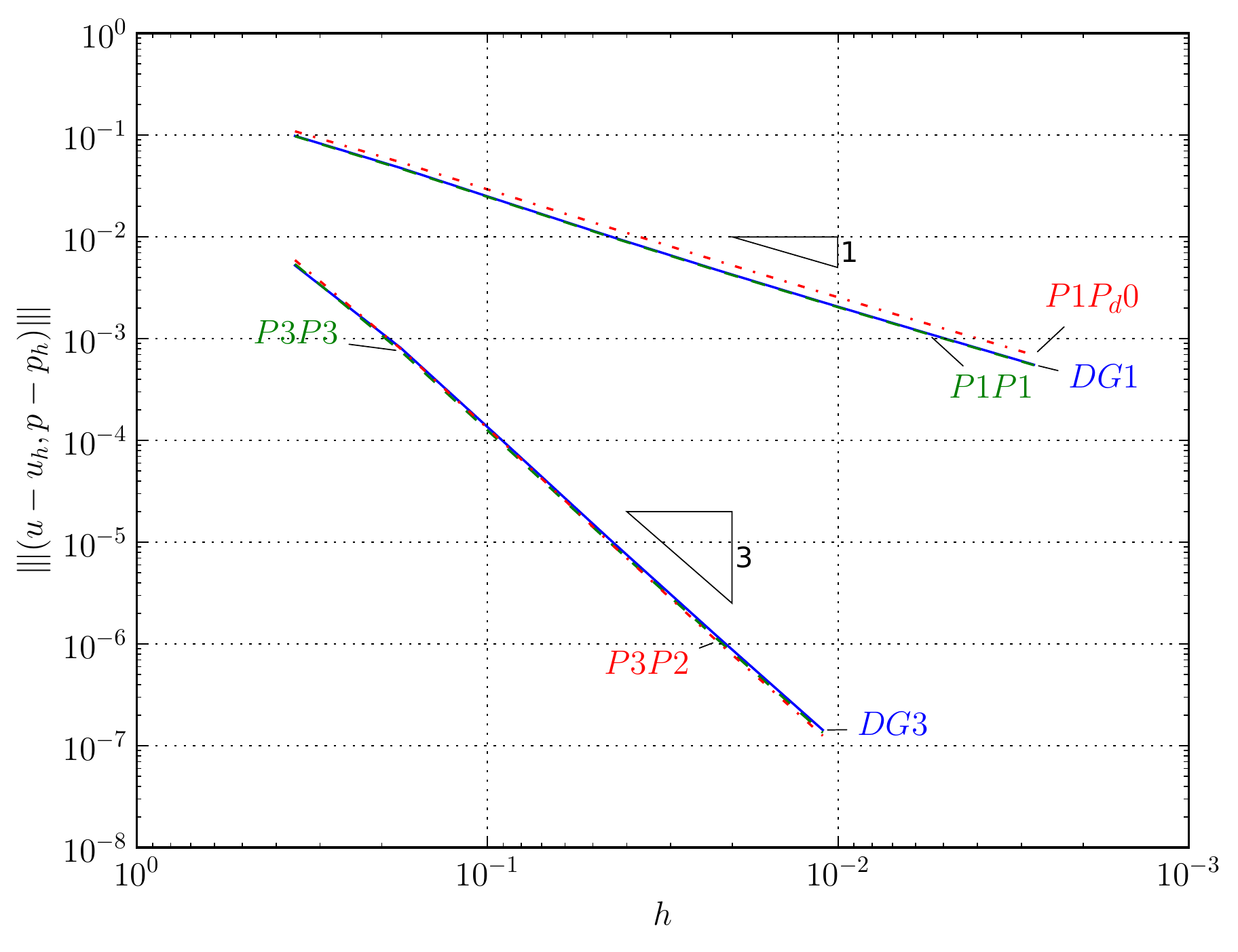}
    \caption{$k=1,3$}
    \label{fig:ref_case1_sto_odd}
  \end{subfigure}
  \begin{subfigure}{0.45\textwidth}
    \centering
    \includegraphics[width=\linewidth]{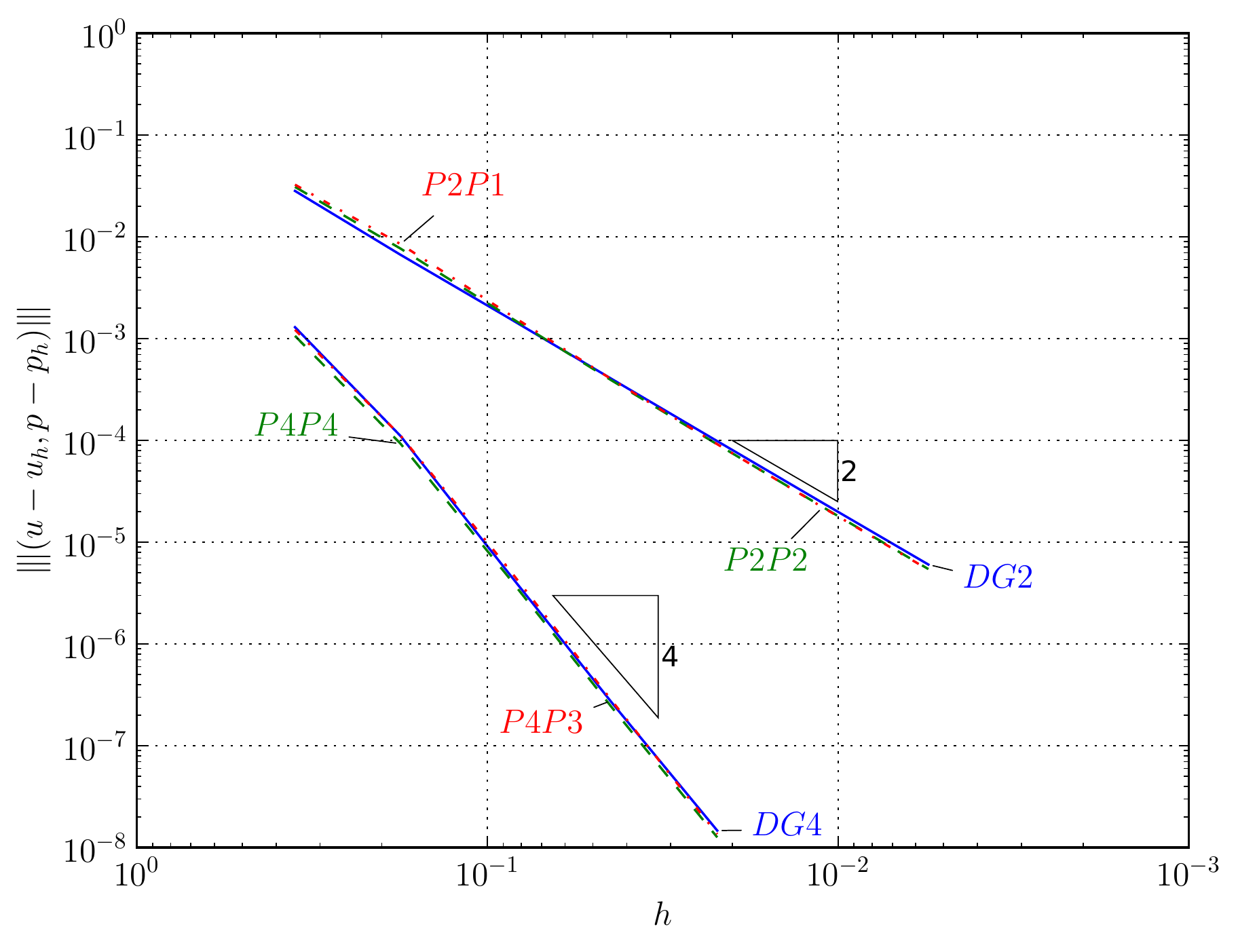}
    \caption{$k=2,4$}
    \label{fig:ref_case1_sto_even}
  \end{subfigure}
  \caption{$\rrrvert (\boldsymbol{u-u_h}, p-p_h) \rrrvert$  for $(\boldsymbol{u},p)$ for $k$ order of $\boldsymbol{u}$ and uniform refinement for a smooth solution.}
  \label{fig:case1_results_sto}
\end{figure}

We start from the mesh in Figure~\ref{fig:ini_mesh1} and perform a uniform refinement of the mesh in each iteration. We plot the error convergence for $\boldsymbol{u}$ and $p$ as a function of the mesh size $h$ in logarithmic scale.  Figures~\ref{fig:case1_results_L2u} and~\ref{fig:case1_results_L2p} show the errors in $L^2$ norm for $\boldsymbol{u}$ and $p$ for different velocity polynomial orders $k$, while Figure~\ref{fig:case1_results_sto} shows the DG error norms. Figure~\ref{fig:case1_results_L2u} displays the different convergence rates for the error in the velocity field for different continuous trial spaces when compared to the $DG$ solution for all orders $k$. Even though the odd cases ($k=1,3$) converge with the same rate as the DG cases, the even cases ($k=2,4$) do not achieve optimality. Their convergence rate in the form $\|\boldsymbol{u-u_h}\|_{L^2}$ is suboptimal and diverges from the $DG2$ and $DG4$ cases. Figure~\ref{fig:case1_results_L2p} exhibits the convergences of $\| p-p_h\|_{L^2}$ for the different $k$ orders. In all cases, the convergence rate is asymptotically consistent with the one displayed by the $DG$ cases. Except for $k=1$, where the convergence rate is 1 for the $P1P0$ case. Nevertheless, it is higher for the $DG1$ and $P1P1$ cases.  Figure~\ref{fig:case1_results_sto} depicts the $\rrrvert(\boldsymbol{u-u_h},p-p_h) \rrrvert$ convergence for the solutions. We illustrate that for this norm, convergence is optimal in all cases.

\begin{rmrk}[Sub-optimal $L^2$ convergence rate for even orders]

  Figure~\ref{fig:case1_results_L2u} shows the method does not guarantee convergence in $\|\boldsymbol{u-u_h}\|_{L^2}$. We use the symmetric interior penalty Galerkin (SIPG) method for the diffusive terms, but the continuous velocity trial space zeroes some terms from~\eqref{a_sip_bilinear}, which result in suboptimal convergence rates for the velocity for even degree polynomials in the $L^2$. Similarly, the incomplete interior penalty Galerkin (IIPG) method and the non-symmetric interior penalty Galerkin (NIPG)~\cite{doi:10.1137/1.9780898717440}.

\end{rmrk}

\begin{figure}[H]
\centering
	\begin{subfigure}{0.45\textwidth}
		\includegraphics[width=\linewidth]{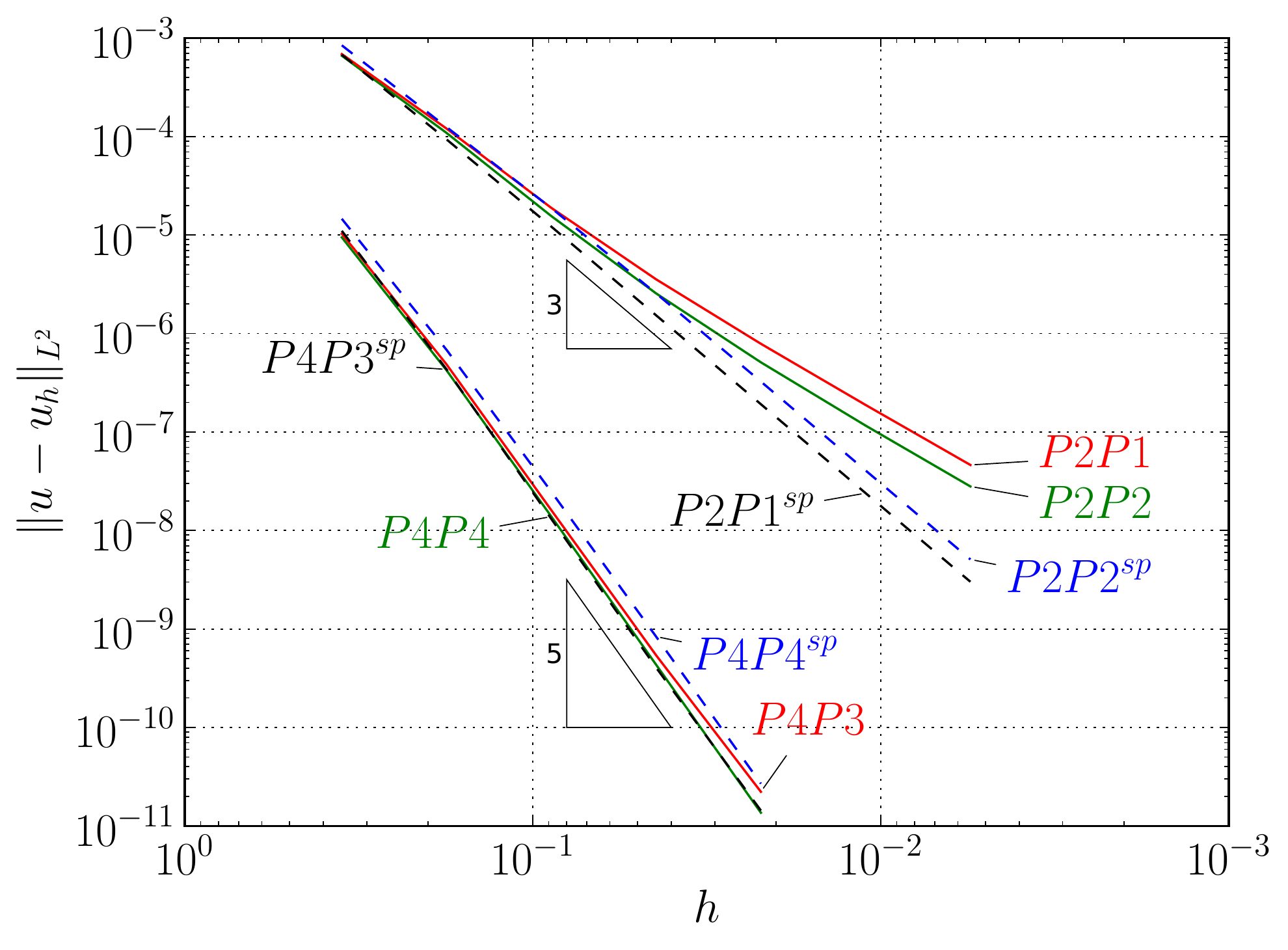}
		\caption{\small $\| \boldsymbol{u-u_h}\|_{L^2}$}
		\label{fig:case1_results_L2u_sp1}
	\end{subfigure}
	\begin{subfigure}{0.45\textwidth}
		\includegraphics[width=\linewidth]{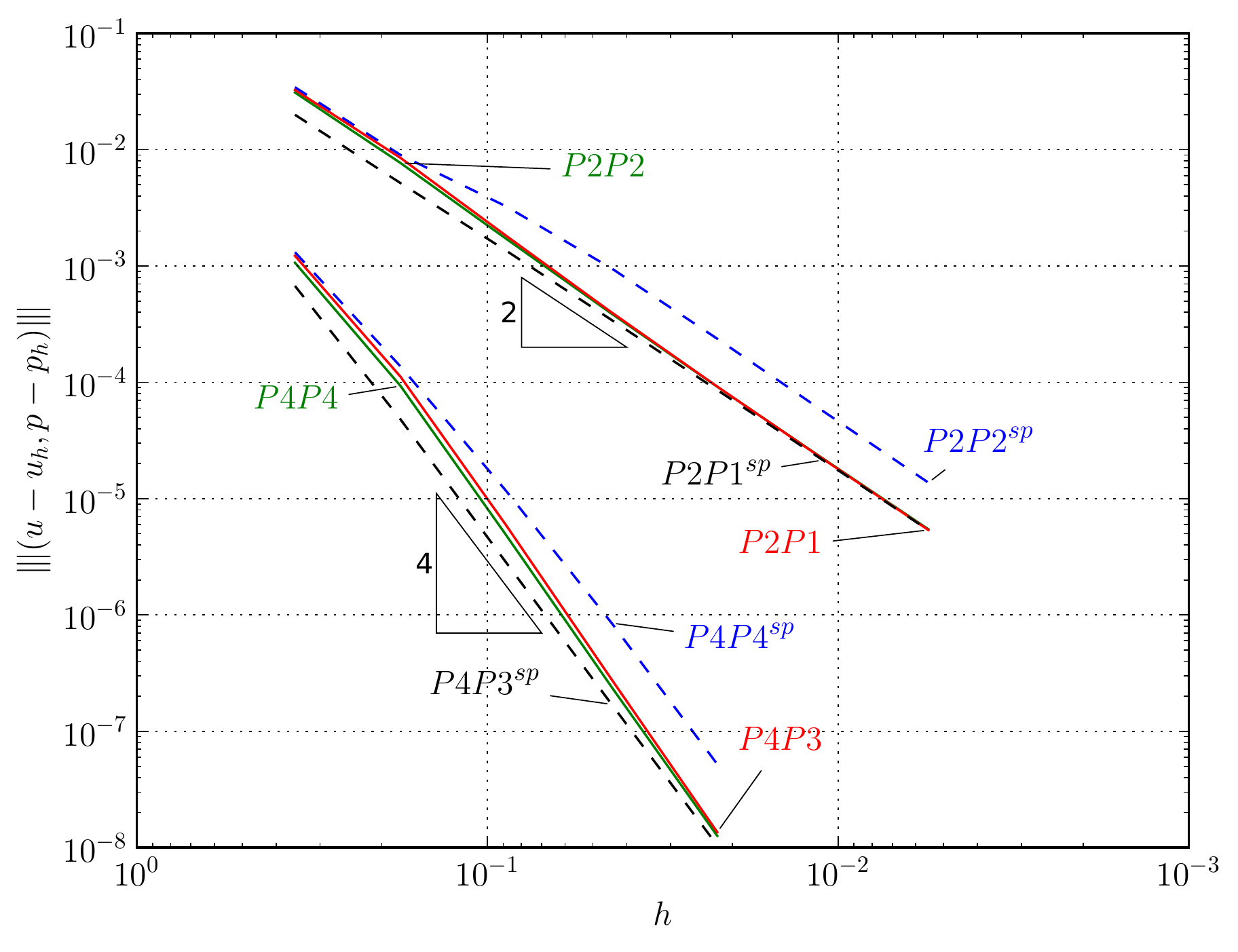}
		\caption{\small $\rrrvert (\boldsymbol{u-u_h}, p-p_h) \rrrvert$ }
		\label{fig:case1_results_sto_sp}
	\end{subfigure}
	
	\caption{\small Error for $\boldsymbol{u}$ and $V_h$ with $k=2,4$ order for $\mathbf{u}$ and $\beta=3$  for the super-penalized case.}
	\label{fig:case1_result_sp}
\end{figure}

\begin{rmrk}[Recovery of optimal $L^2$ convergence rate for even
  orders via super-penalization]

  We super-penalize the velocity jumps in the residual representation to recover optimal convergence, as proposed for interior penalty methods~\cite{ doi:10.1137/1.9780898717440}. That is, in~\eqref{a_sip_bilinear} we substitute
$$
+ \sum_{F \in F_h} \frac{\eta}{h_F} \int_{F} \llbracket \boldsymbol{v_{h}} \rrbracket \cdot \llbracket \boldsymbol{w_{h}} \rrbracket 
$$
by
$$
+ \sum_{F \in F_h} \frac{\eta}{h^{\beta}_F} \int_{F} \llbracket \boldsymbol{v_{h}} \rrbracket \cdot \llbracket \boldsymbol{w_{h}} \rrbracket 
$$
where $\beta \geq 1$ but to set its, we follow the condition $\beta(d-1) \geq 3$ from the interior penalty methods~\cite{ doi:10.1137/1.9780898717440}. Thus, we solve again the even cases ($k=2,4$) using super-penalization and set $\beta=3$.  Figure~\ref{fig:case1_result_sp} compares the convergence of the different cases with and without super-penalization. We label the cases with $\beta=3$ as \textit{"sp"} where setting $\beta$ to 3 restores optimality for the velocity field while preserving optimal convergence in the DG norm. Nevertheless, we do not recommend super-penalization as the condition number of the algebraic system suffers significantly hampering the ability of the iterative solver to converge. More importantly, the adaptive scheme restores optimal convergence without requiring an additional tuning parameter.

\end{rmrk}

\begin{figure}[h]
  \centering
  \begin{subfigure}{0.495\textwidth}
    \centering
    \includegraphics[width=1.15\linewidth]{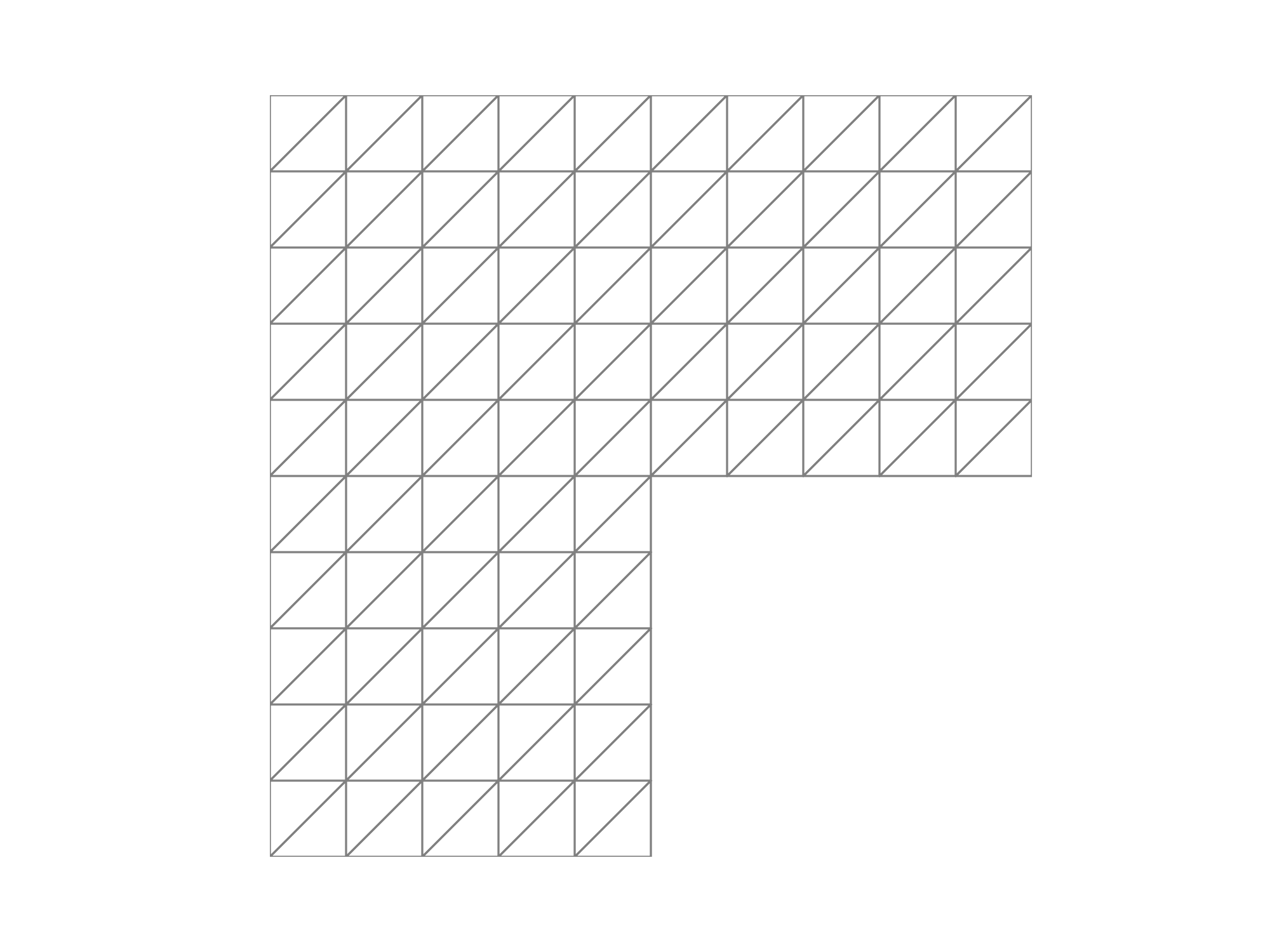}
    \caption{Level 0}
    \label{fig:case2_mesh0}
  \end{subfigure}
  \begin{subfigure}{0.495\textwidth}
    \centering
    \includegraphics[width=1.15\linewidth]{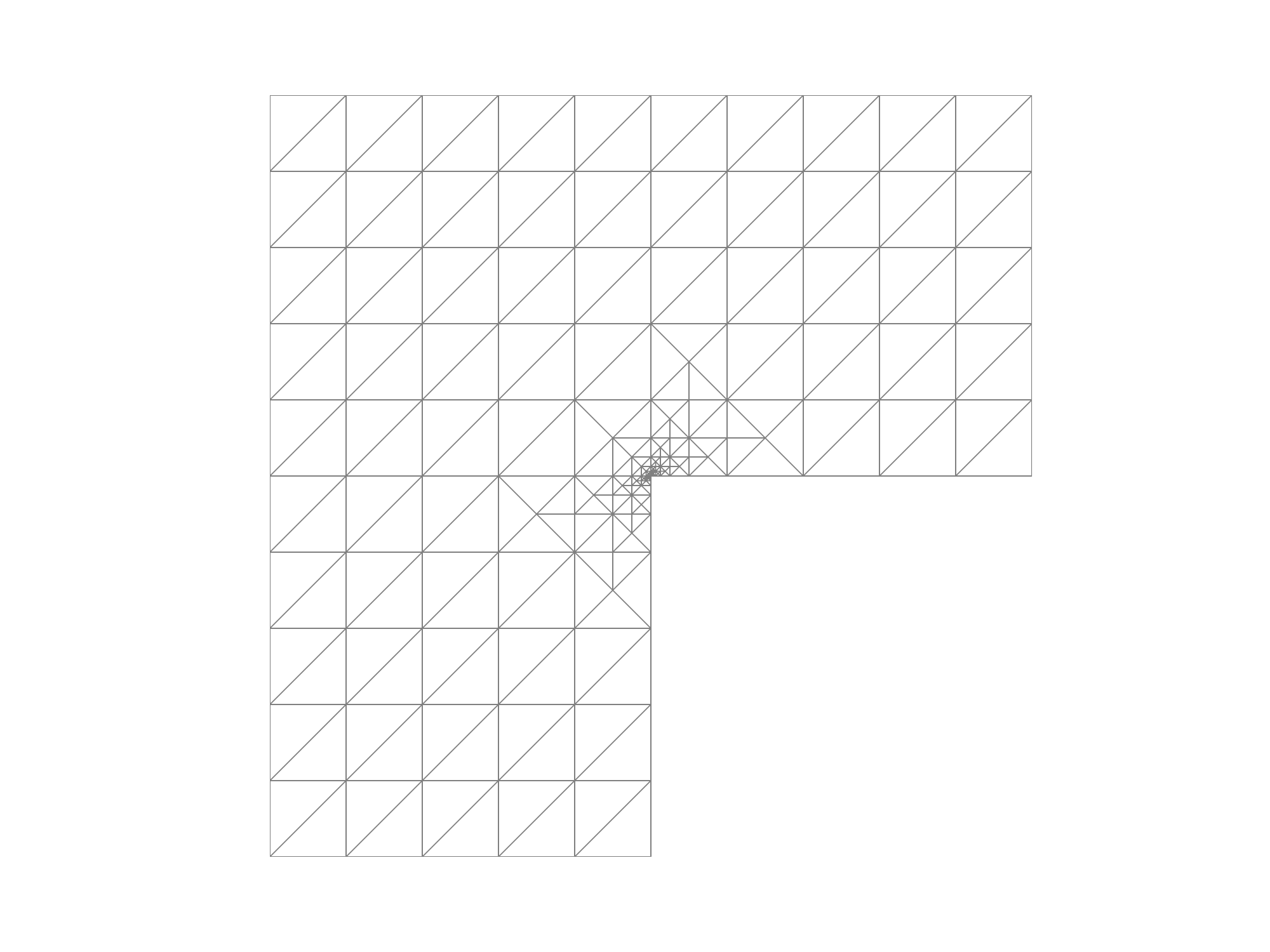}
    \caption{Level 10}
    \label{fig:case2_mesh10}
  \end{subfigure} 
  \begin{subfigure}{0.495\textwidth}
    \centering
    \includegraphics[width=1.15\linewidth]{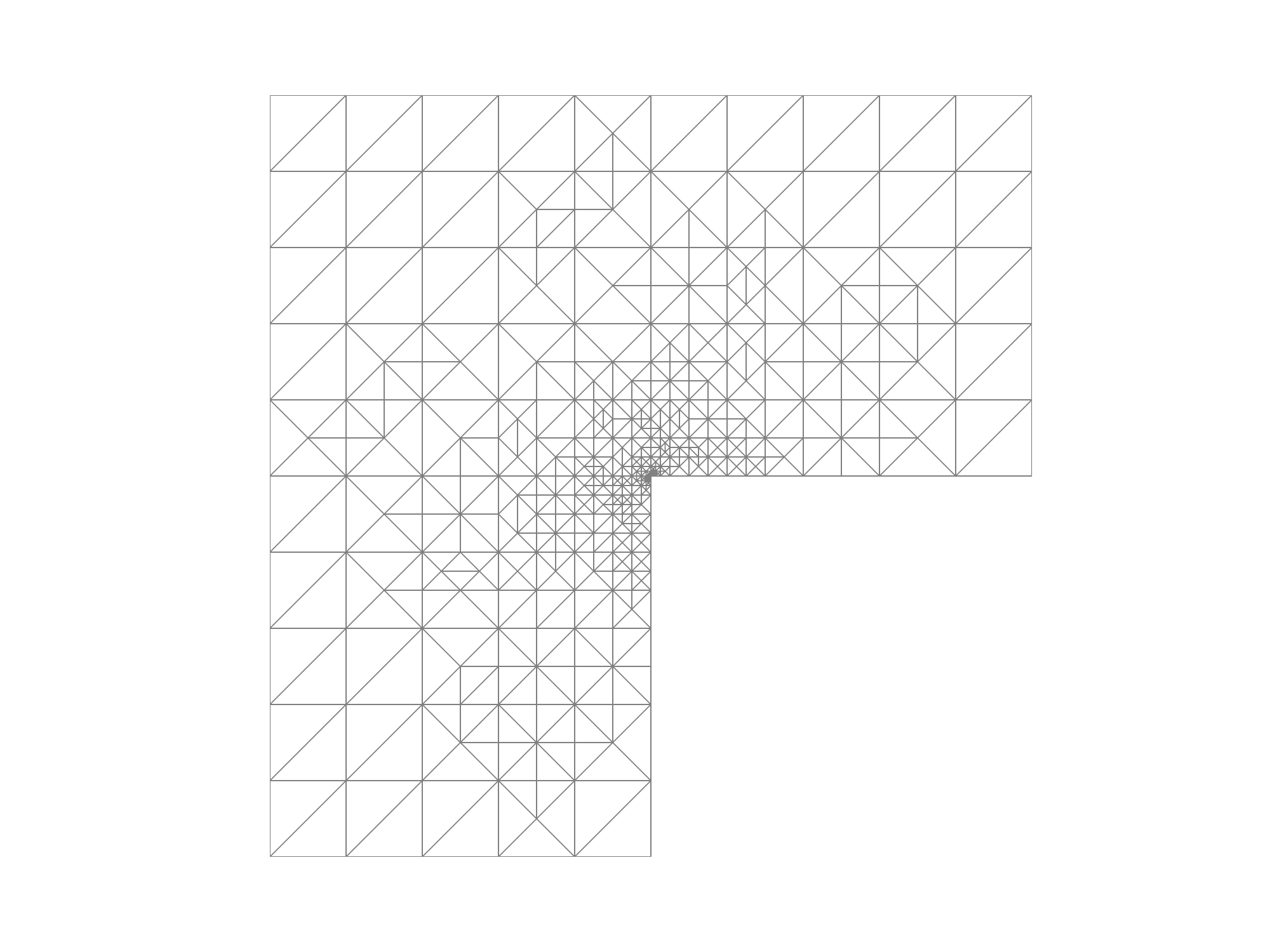}
    \caption{Level 20}
    \label{fig:case2_mesh20}
  \end{subfigure}
  \begin{subfigure}{0.495\textwidth}
    \centering
    \includegraphics[width=1.15\linewidth]{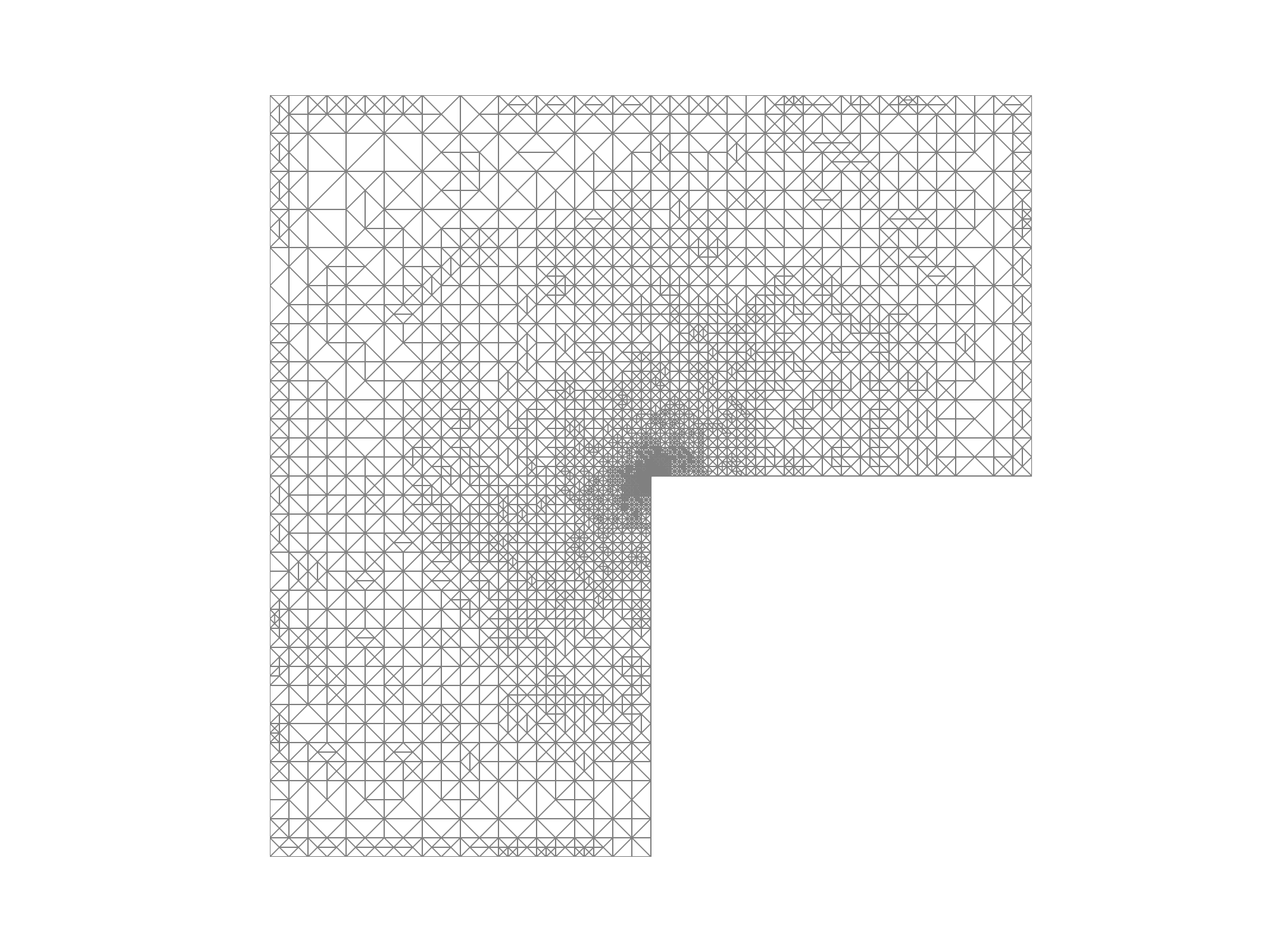}
    \caption{Level 20}
    \label{fig:case2_mesh40}
  \end{subfigure}
  
  \caption{$P1P1$ adaptive refinement solution for Case 2. }
  \label{fig:case2_evp1_mesh}
\end{figure}

\subsection{L-shape with adaptive refinement}

For the adaptive mesh refinement case, we consider the L-shape $\Omega = (-1,1)^2 \backslash ((0,1)\times (-1,0)) $ with analytical solution for $\boldsymbol{u}$ and $p$ ~\cite{10.2307/2153481} as
%
%
\begin{equation}
\boldsymbol{u}_{sol} = 
\begin{bmatrix}
r^{\alpha}[(1+\alpha)\sin\varphi \psi(\varphi)+\cos(\varphi)\partial_{\varphi}\psi(\varphi)],
\\
r^{\alpha}[\sin(\varphi)\partial_{\varphi}\psi(\varphi)- (1+\alpha)\cos\varphi \psi(\varphi)]

\end{bmatrix},   
\end{equation}

\begin{equation}\label{p_exact_2}
p_{sol} = -r^{\alpha-1}[(1+\alpha)^{2}\partial_{\varphi}\psi(\varphi)+\partial^3_{\varphi}\psi(\varphi)]/(1-\alpha)
\end{equation}
with
\begin{equation}
	\begin{split}
		\psi(\varphi) = &\sin(( 1 + \alpha)\varphi)\cos(\alpha\omega)/( 1-\alpha) - \cos(( 1 +\alpha)\varphi)\\&+ sin((\alpha - 1)\varphi)\cos(\alpha\omega)/(1-\alpha) + cos((\alpha-1)\varphi),
		\\& \alpha= 856399/1572864,\qquad \omega=\frac{3\pi}{2}.
	\end{split}
\end{equation}

We set the forcing term equal to zero and impose Dirichlet boundary conditions on the entire boundary, setting homogeneous values on the edges of the reentrant corner and nonhomogeneous ones on the complementary part. We consider for this problem the cases $PkPk$ with $k=1,2$ where the test spaces are $V_h \times Q_h = [\Pol^{k}_{d}(\Omega_h)]^2 \times \Pol^{k}_{d,0}(\Omega_h)$ and the trial spaces are $U_h\times P_h= [\mathbb{P}^k(\Omega_h)]^2 \times \mathbb{P}_{0}^k(\Omega_h)$. For the adaptive refinement scheme, we choose a user-specified fraction of $"0.1"$.

\begin{figure}[H]
	\begin{subfigure}{0.5\textwidth}
		\centering
		\includegraphics[width=0.9\linewidth]{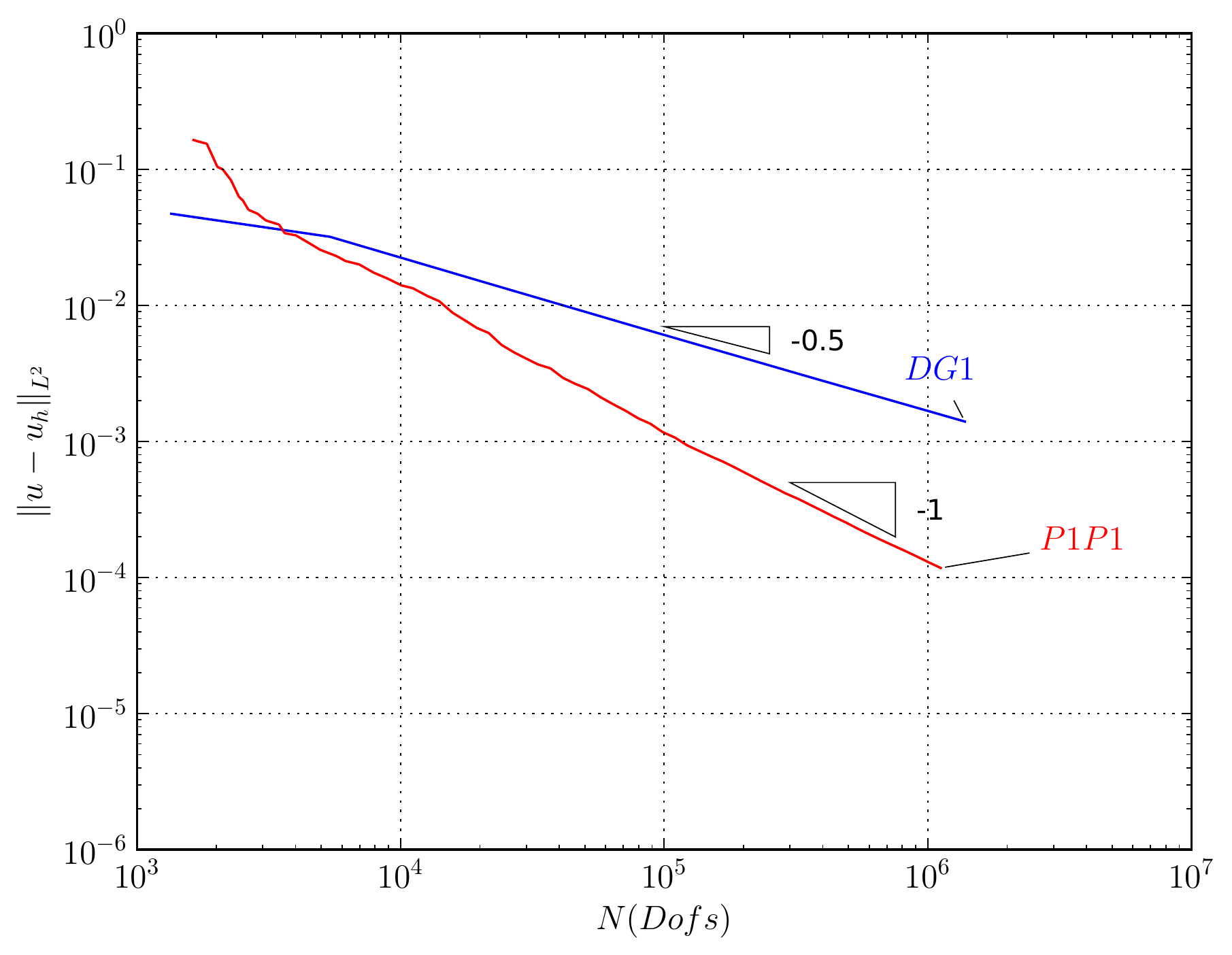}
		\caption{$k=1$}
		\label{fig:ref_case2_u_1}
	\end{subfigure}
	\begin{subfigure}{0.5\textwidth}
		\centering
		\includegraphics[width=0.9\linewidth]{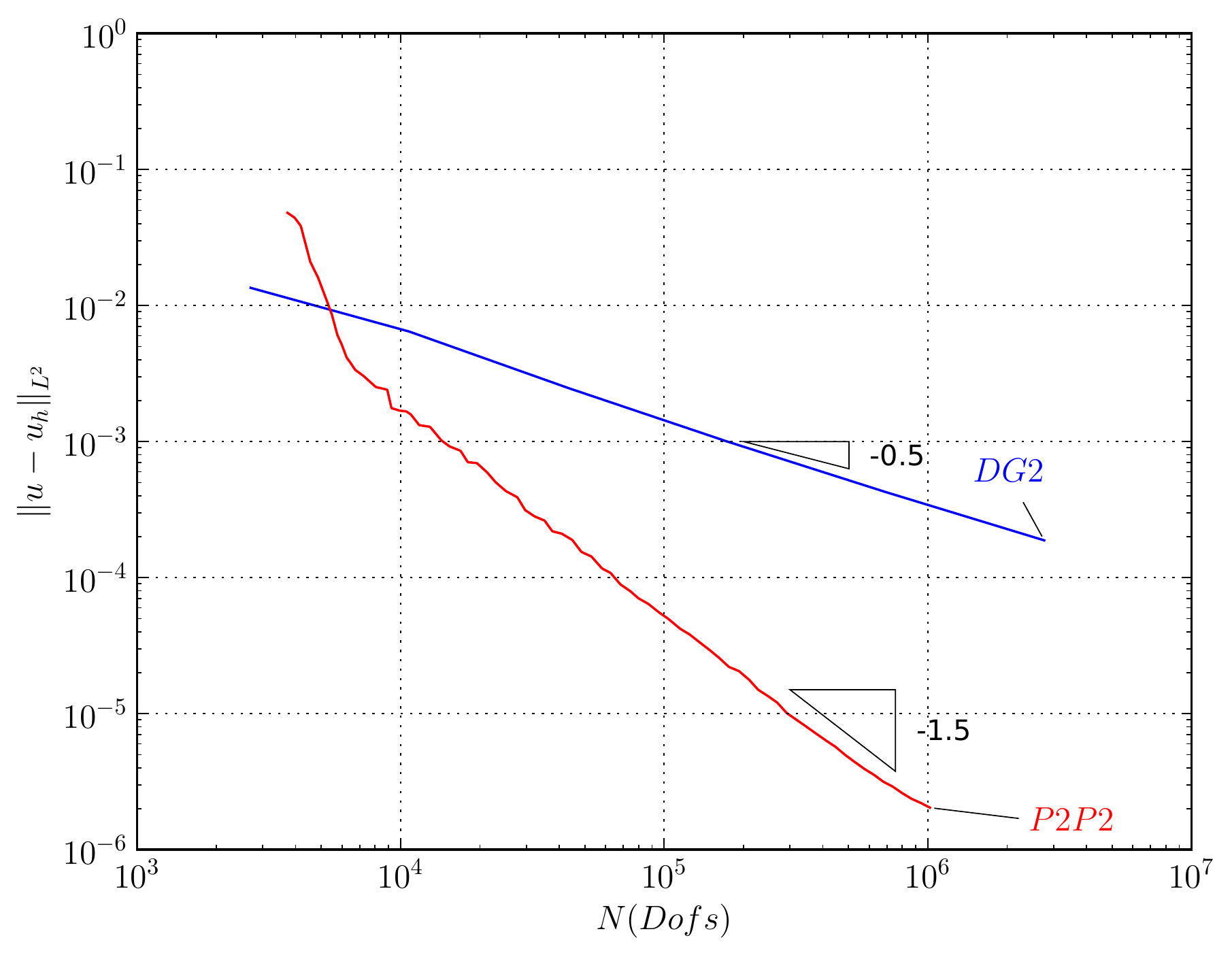}
		\caption{$k=2$}
		\label{fig:ref_case2_u_2}
	\end{subfigure}
	\caption{$\| \boldsymbol{u-u_h}\|_{L^2}$ for $k$ order of $\boldsymbol{u}$ with uniform refinement for $DGk$ and adaptive refinement for $PkPk$ cases.}
	\label{fig:case2_results_L2u}
\end{figure}
\begin{figure}[H]
	\begin{subfigure}{0.5\textwidth}
		\centering
		\includegraphics[width=0.9\linewidth]{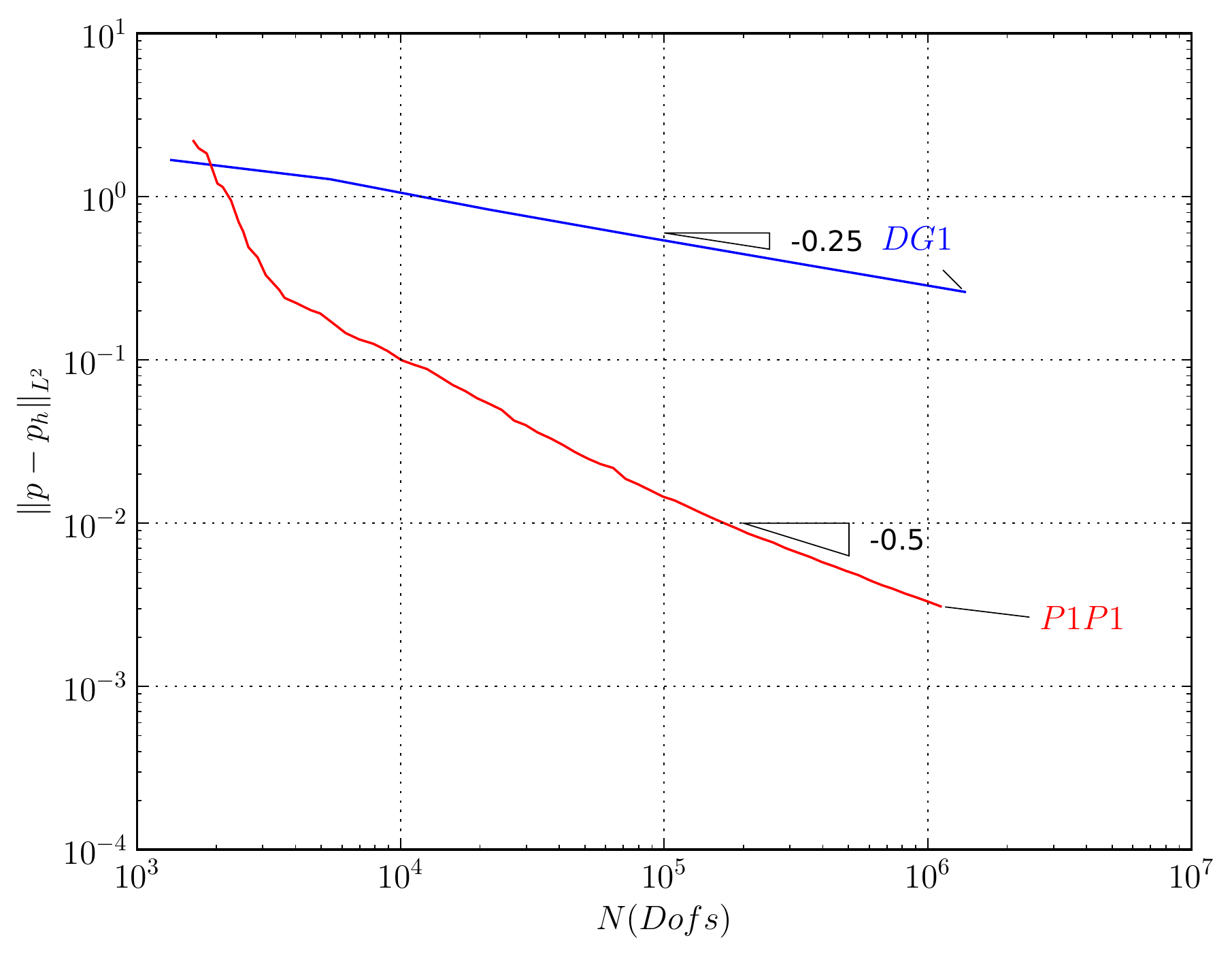}
		\caption{$k=1$}
		\label{fig:ref_case2_p_1}
	\end{subfigure}
	\begin{subfigure}{0.5\textwidth}
		\centering
		\includegraphics[width=0.9\linewidth]{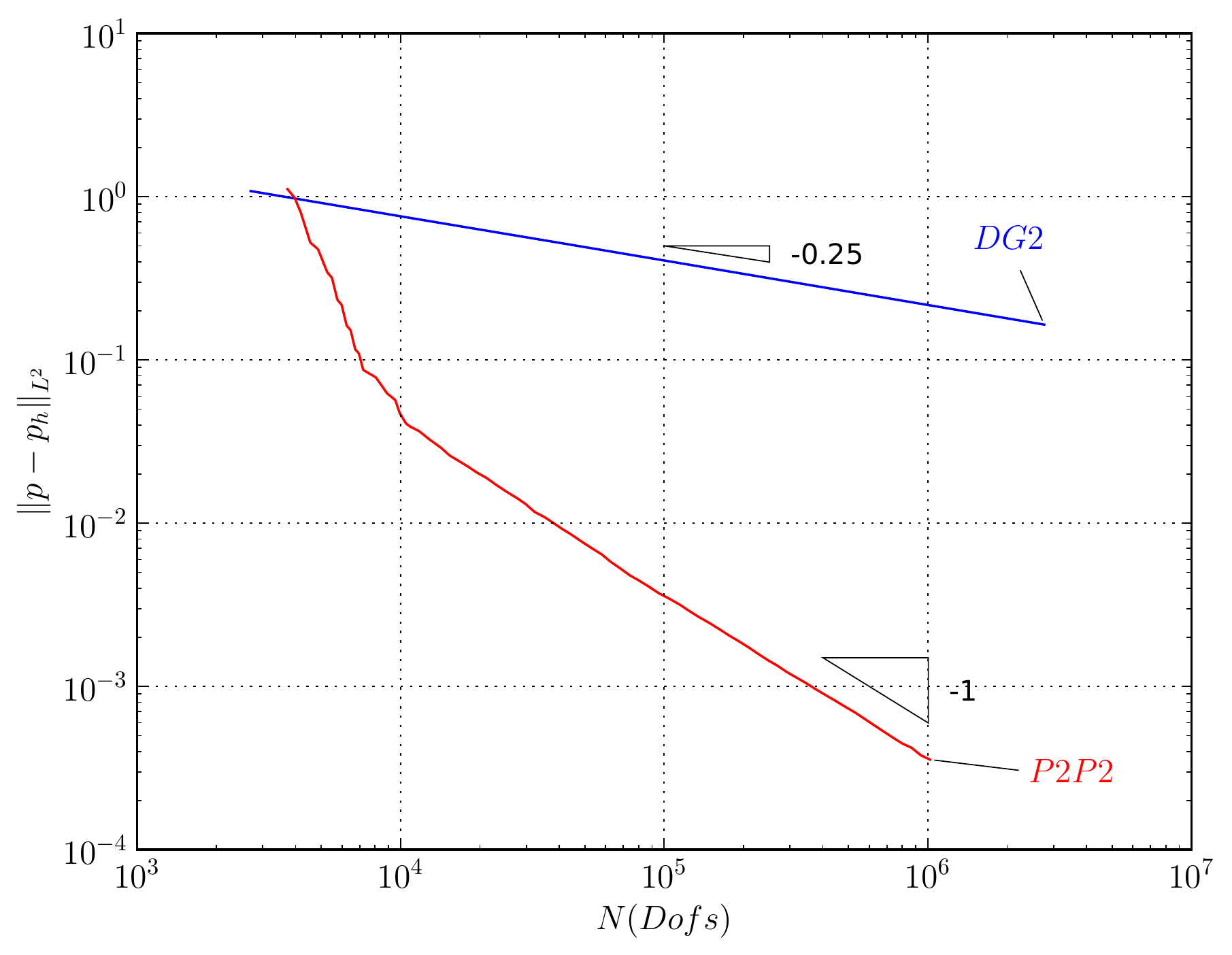}
		\caption{$k=2$}
		\label{fig:ref_case2_p_2}
	\end{subfigure}
	\caption{$\|p-p_h\|_{L^2}$ for $k$ order of $\boldsymbol{u}$ with uniform refinement for $DGk$ and adaptive refinement for $PkPk$ cases.}
	\label{fig:case2_results_L2p}
\end{figure}
\begin{figure}[H]
	\begin{subfigure}{0.5\textwidth}
	\centering
	\includegraphics[width=0.9\linewidth]{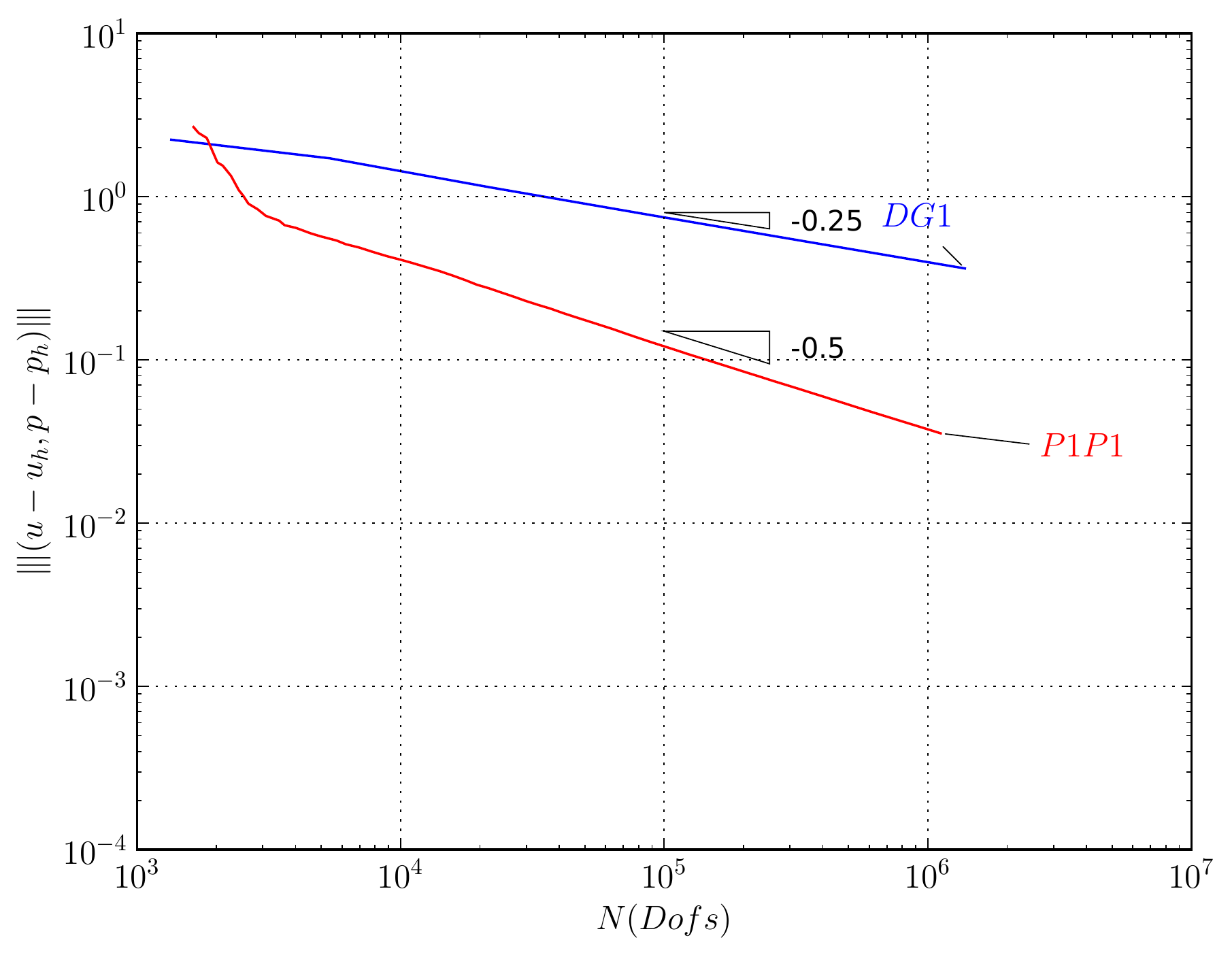}
	\caption{$k=1$}
	\label{fig:ref_case2_sto_1}
\end{subfigure}
\begin{subfigure}{0.5\textwidth}
	\centering
	\includegraphics[width=0.9\linewidth]{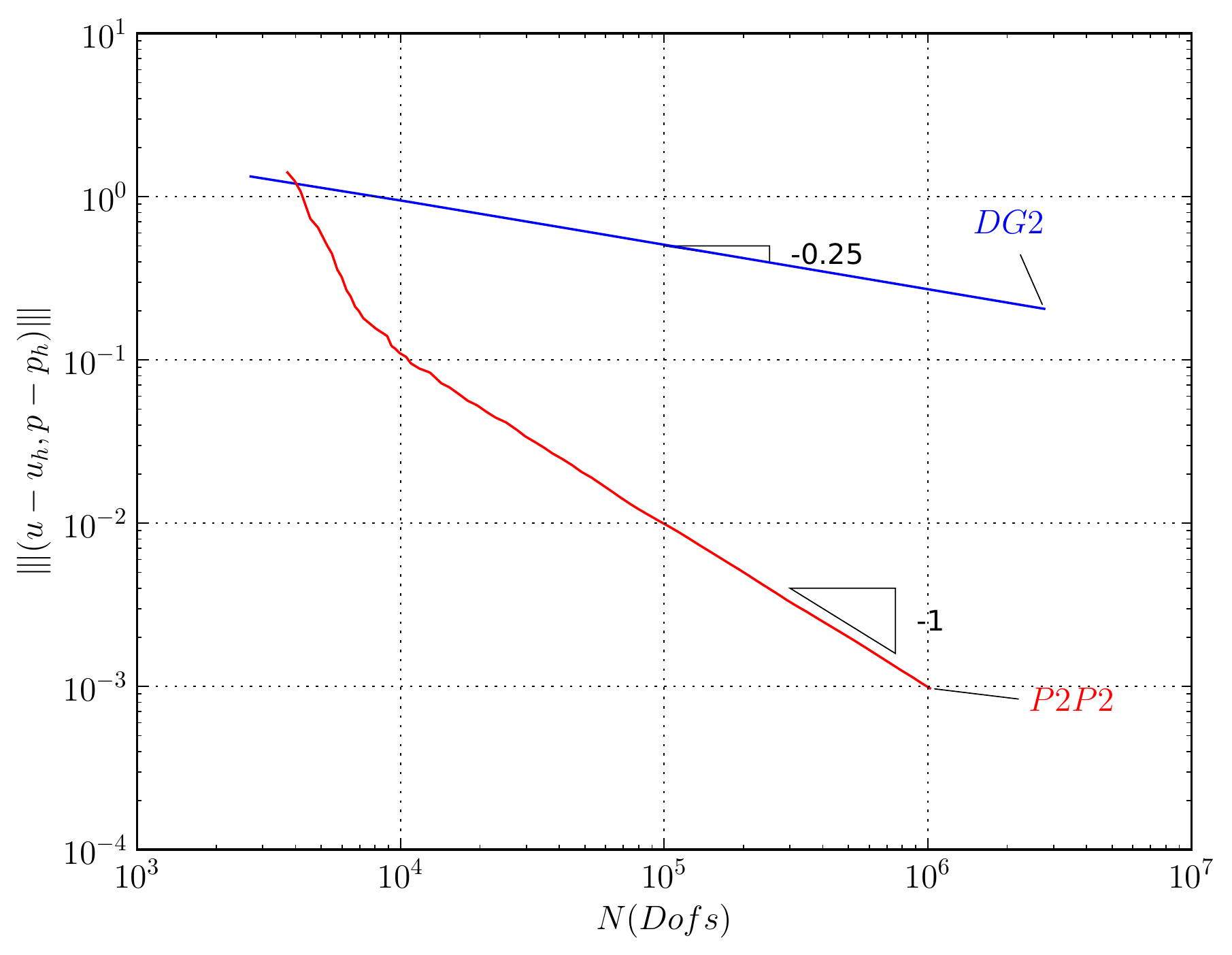}
	\caption{$k=2$}
	\label{fig:ref_case2_sto_2}
\end{subfigure}
	\caption{$\rrrvert (\boldsymbol{u-u_h}, p-p_h) \rrrvert$  for $(\boldsymbol{u},p)$ for $k$ order of $\boldsymbol{u}$ with uniform refinement for $DGk$ and adaptive refinement for $PkPk$ cases.}
	\label{fig:case2_results_sto}
\end{figure}

We start the refinement process from a coarse, uniform mesh, see Figure~\ref{fig:case2_mesh0}, and apply an adaptive refinement driven by the error representations we introduce. We compare these results against the discontinuous Galerkin solutions $DG1$ and $DG2$ with uniform refinement. Figure~\ref{fig:case2_evp1_mesh} shows the meshes for different levels for different refinement levels for the $P1P1$ case.  Figure~\ref{fig:case2_results_L2u} shows $\| \boldsymbol{u-u_h}\|_{L^2}$ versus the numbers of $DOFs$ of the solution. Automatic adaptive refinement delivers a clear improvement in the convergence for the $PkPk$ cases against the $DGk$ cases. Figures~\ref{fig:case2_results_L2p} and~\ref{fig:case2_results_sto} show $\| \boldsymbol{p-p_h}\|_{L^2}$ and $\rrrvert (\boldsymbol{u-u_h}, p-p_h) \rrrvert$, respectively, versus the numbers of $DOFs$ of the solution. In both error norms, the $PkPk$ adaptive cases converge faster to the solution than the standard $DGk$ cases with uniform refinement.

\begin{figure}[h!]
  \centering
  \includegraphics[width=0.6\linewidth]{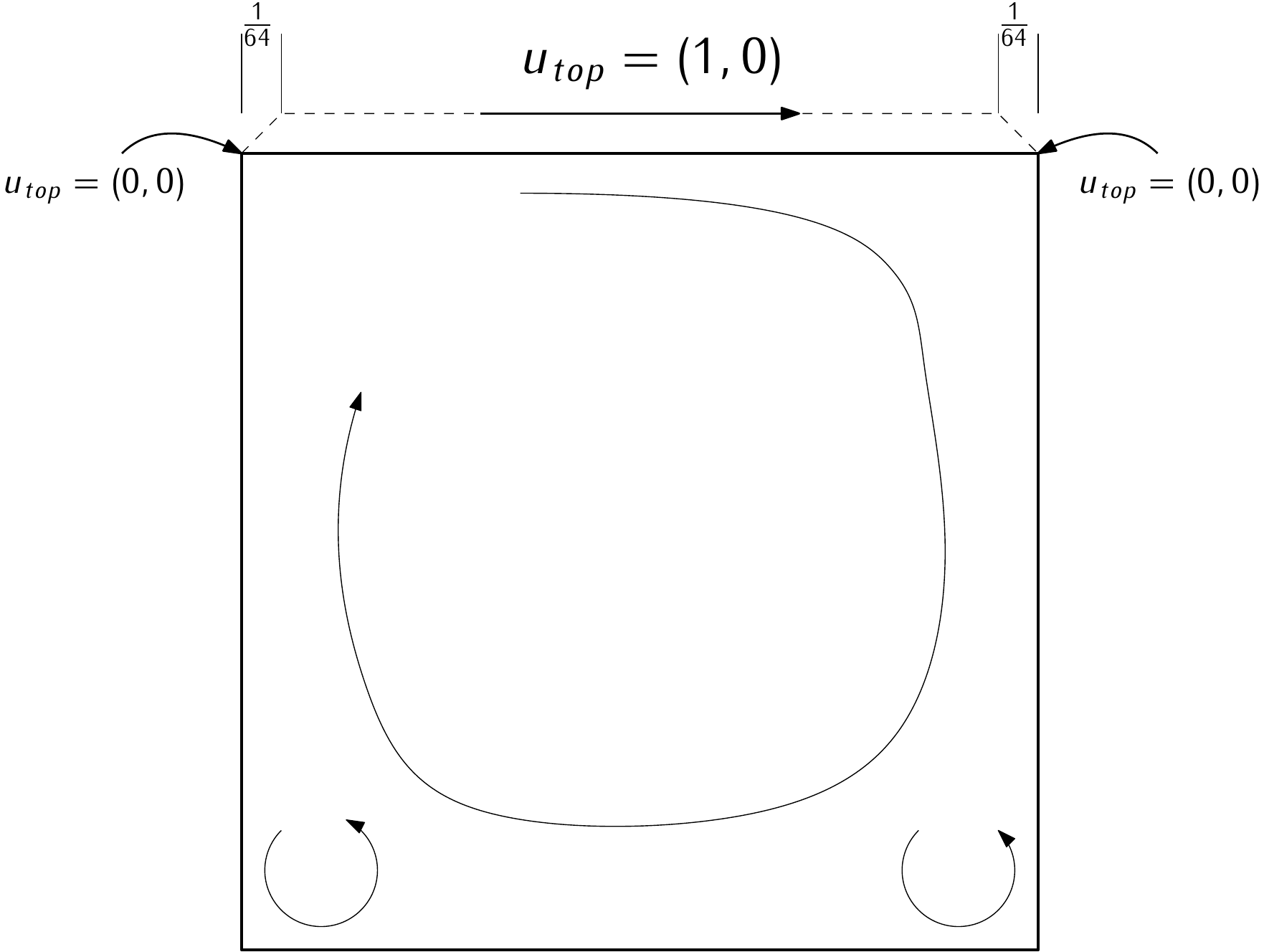}
  \caption{Boundary conditions for the lid-cavity case in 2D.}
  \label{fig:lidCavity}
\end{figure}

\subsection{Lid-driven cavity flow}

As a last example, we consider the well-known lid-driven cavity flow
problem (see~\cite{LIU2006580,CORTES2015123} ). We set the source term
to $\boldsymbol{f}=0$ and consider no-slip boundary conditions on the
bottom, left, and right boundaries ($\boldsymbol{u} = (0,0)$). At the
top, as Figure~\ref{fig:lidCavity} shows, we impose the velocity profile $\boldsymbol{u}_{top}$~(see~\cite{ ROBERTS20172018}). Therefore, we establish a linear transition from 0 to 1 with a slope  $\epsilon = 1/64$ in the top corners to avoid a discontinuity.

\begin{figure}
    \centering
	\begin{subfigure}{0.48\textwidth}
		\centering
		\includegraphics[height=\textwidth]{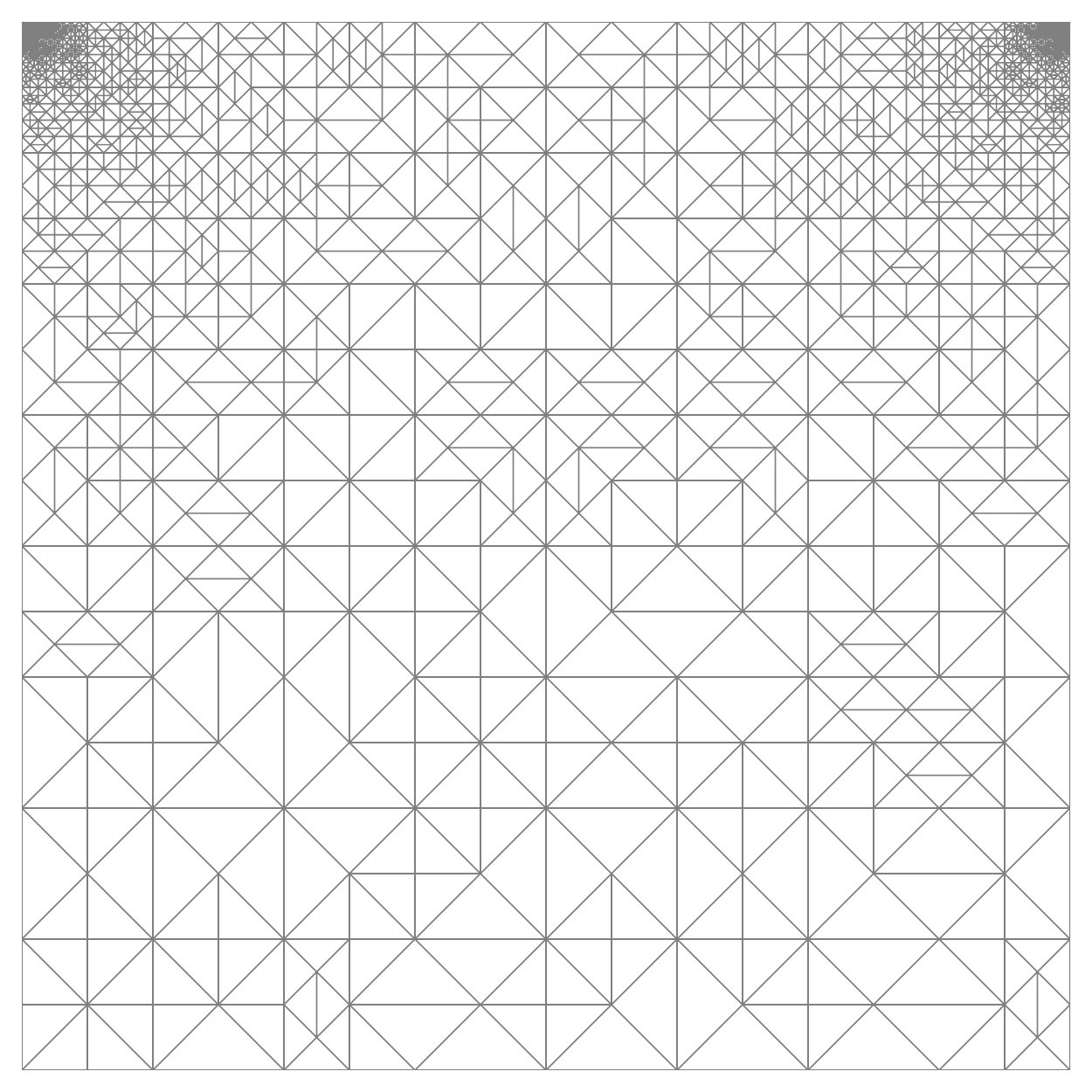}
		\caption{Mesh}
		\label{fig:lid_mesh}
	\end{subfigure}
	\begin{subfigure}{0.48\textwidth}
		\centering
	\includegraphics[height=\textwidth]{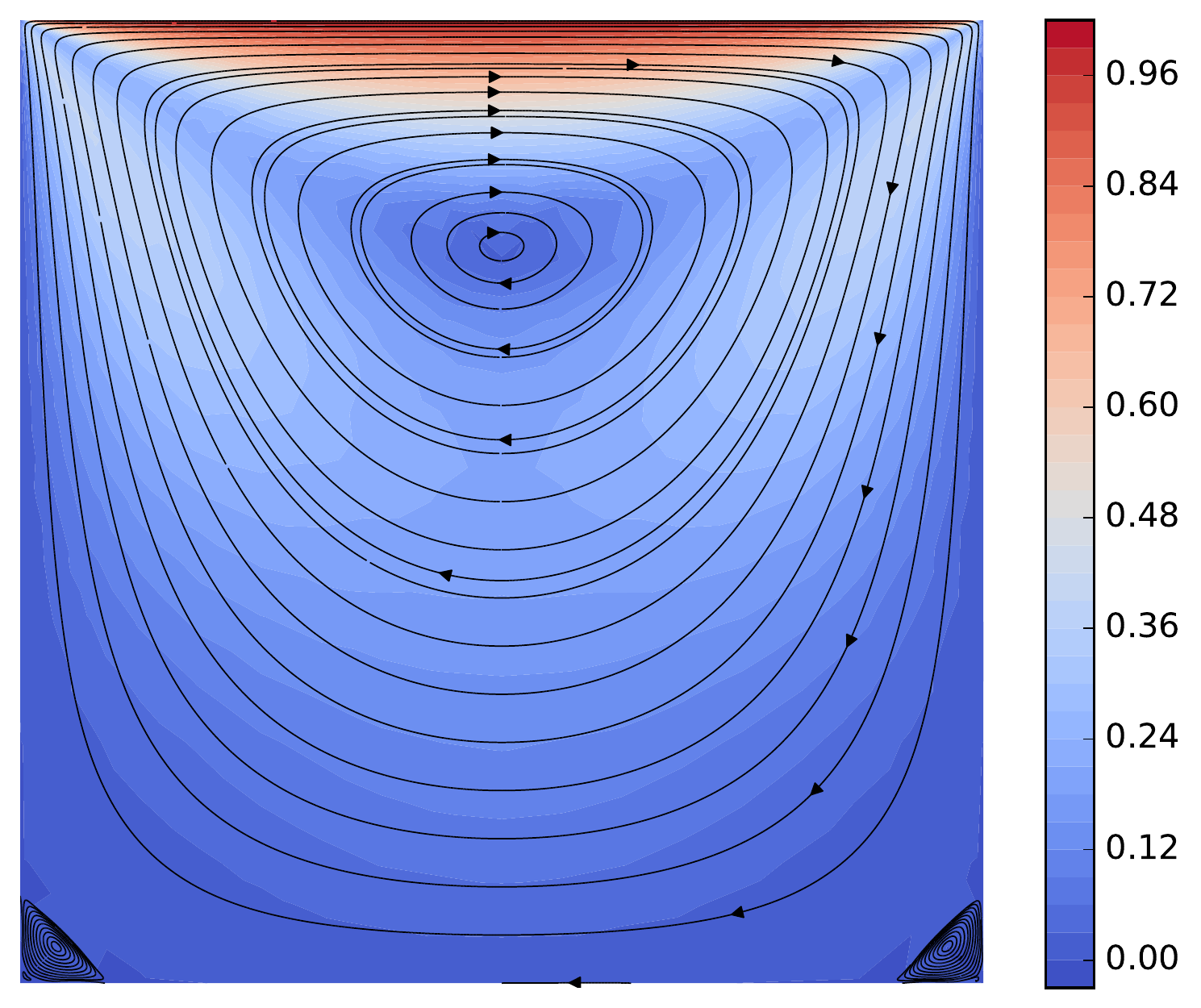}
	\caption{$\mathbf{u}$ solution}
	\label{fig:lid_u}
\end{subfigure} \\
	\begin{subfigure}{0.8\textwidth}
		\centering
		\includegraphics[width=\linewidth]{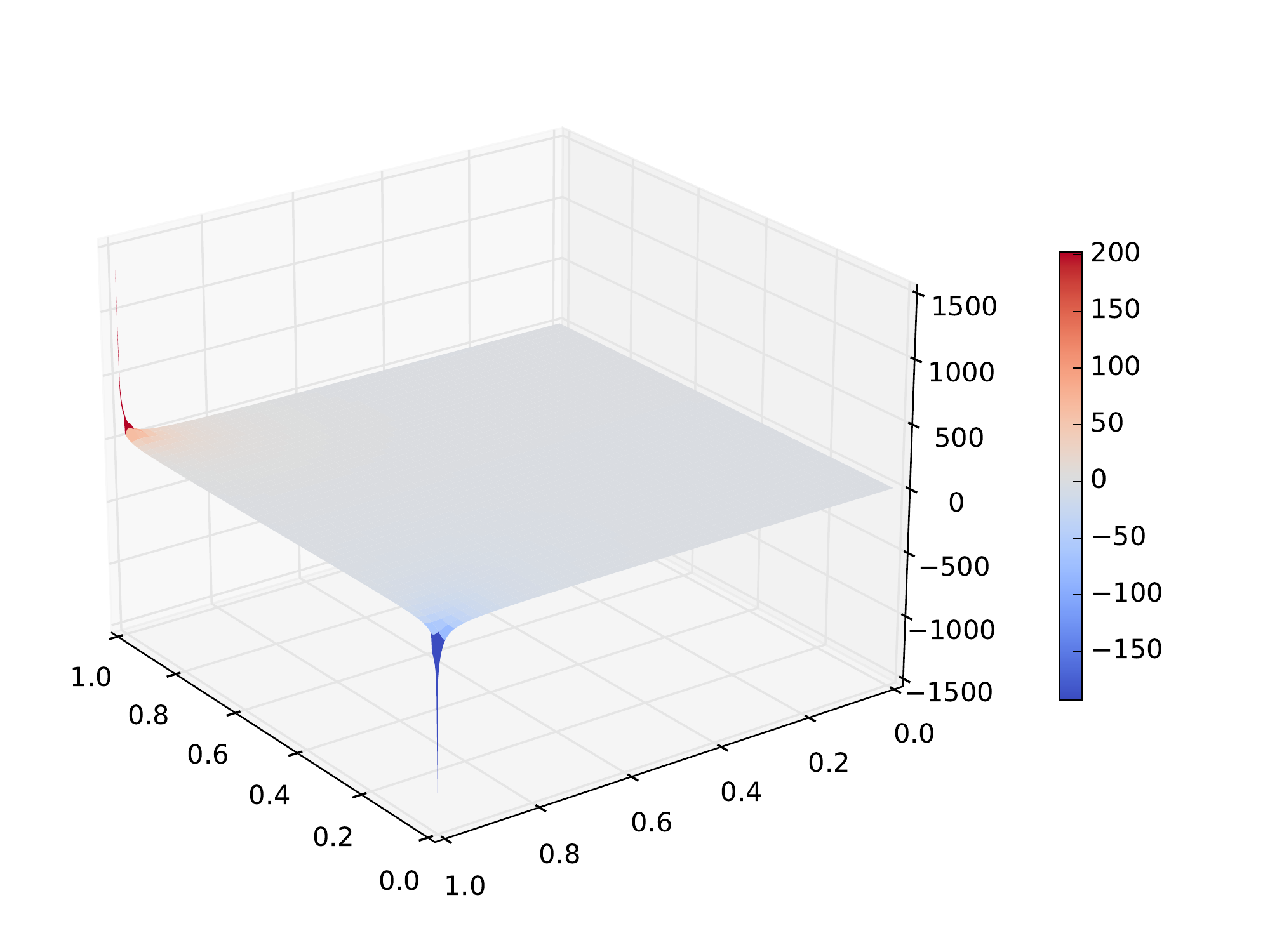}
		\caption{$p$ solution}
		\label{fig:lid_p}
	\end{subfigure}
	
	\caption{$P3P3$ automatically adapted solution for the lid-driven cavity,.  $\boldsymbol{u}$ and $p$ fields correspond to refinementlevel 50.}
	\label{fig:lid_p2p1_mesh}
\end{figure}

We use test spaces $V_h \times Q_h = [\Pol^{3}_{d}(\Omega_h)]^2 \times \Pol^{3}_{d,0}(\Omega_h)$ and the trial spaces $U_h\times P_h= [\mathbb{P}^3(\Omega_h)]^2 \times \mathbb{P}_{0}^3(\Omega_h)$. 
Figure~\ref{fig:lid_p2p1_mesh} presents the solution for $\mathbf{u}$ and $p$ at the refinement level 50. 

\begin{figure}
		\centering
		\includegraphics[height=0.5\textwidth]{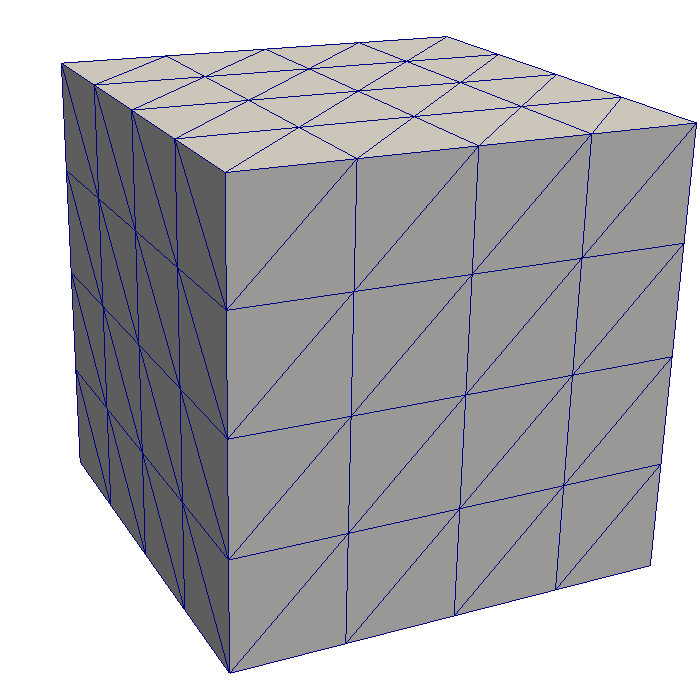}	
	\caption{Initial mesh for the lid-driven cavity 3D case.}
	\label{fig:lid_3D_mesh}
\end{figure}

\begin{figure}
	\centering
	\begin{subfigure}{0.48\textwidth}
		\centering
		\includegraphics[height=0.8\textwidth]{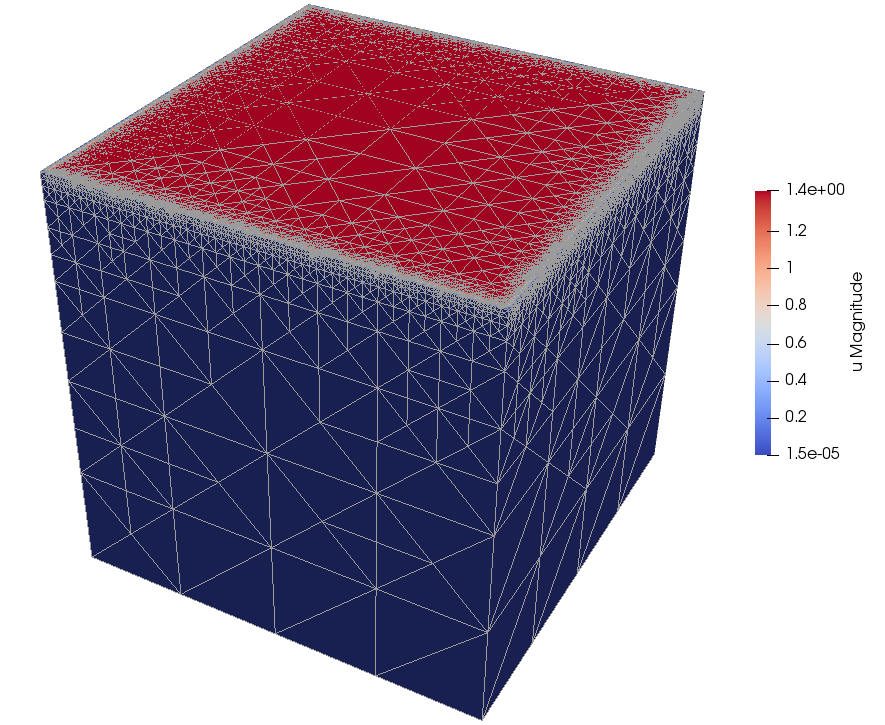}
		\caption{The $\mathbf{u}$ solution and mesh.}
		\label{fig:lid_3D_u}
	\end{subfigure}
	\begin{subfigure}{0.48\textwidth}
		\centering
		\includegraphics[height=0.8\textwidth]{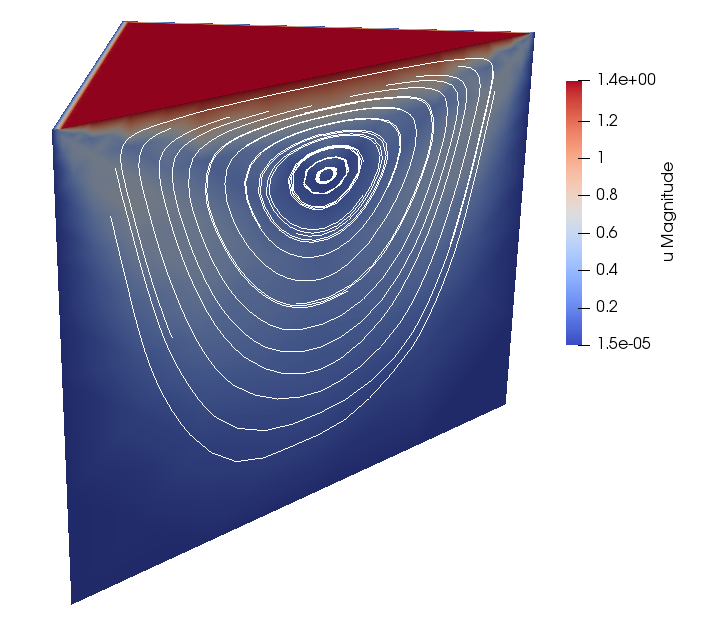}
		\caption{The $\mathbf{u}$ solution and streamlines in the middle of the domain.}
		\label{fig:lid_3D_u2}
	\end{subfigure} 
	
	\caption{$P3P3$ automatically adapted solution for the 3D lid-driven cavity for the refinement level 26.}
	\label{fig:lid_3D}
\end{figure}

Furthermore, we consider the 3D lid-driven cavity flow. For this case, we consider the domain $\Omega = (0,1)^3$ and consider the following boundary condition for the top face: $\mathbf{u}_{top} = (1,1,0)$ with the same linear variation near the edges like the 2D case. We set homogeneous boundary conditions for the rest of the domain's faces. Figure~\ref{fig:lid_3D_mesh} shows the initial mesh while Figure~\ref{fig:lid_3D} presents the solution for the refinement level 26. 

\section{Conclusions}\label{conclusions}

We use an adaptive-stabilized finite element method (AS-FEM) to develop robust discretizations for the Stokes problem using standard FEM triangular elements as our continuous solution space. We minimize the residual on a dual discontinuous Galerkin (DG) norm and seek the approximate solution fields as continuous subspaces of the DG ones. As a result, we obtain a robust discretization that delivers optimal convergence rates on uniform meshes, except for the $L^2$ norm for the velocity where the polynomial degree is even. For those cases, we show that it is possible to recover optimality by applying super-penalization methods, but we do not recommend their usage. Furthermore, as a result of using  AS-FEM, the user can independently minimize the residual representation of the velocity or the pressure fields in the DG space to perform adaptive refinements of the solutions. 

\section*{Acknowledgements}

This publication was made possible in part by the CSIRO Professorial Chair in Computational Geoscience at the Faculty of Science and Engineering at Curtin University and the Deep Earth Imaging Enterprise Future Science Platforms of the Commonwealth Scientific Industrial Research Organisation, CSIRO, of Australia. This project has received funding from the European Union's Horizon 2020 research and innovation programme under the Marie Sklodowska-Curie grant agreement No 777778 (MATHROCKS). The Curtin Corrosion Centre and the Curtin Institute for Computation kindly provide ongoing support.

\bibliography{mybibfile}

\end{document}